\documentclass[12pt]{article}
\usepackage[english]{babel}
\usepackage{geometry} 
\usepackage{color}

\usepackage{latexsym} % For \Box to end proofs.
\usepackage{amsmath,amssymb,amsthm}
\usepackage[applemac]{inputenc}
 \newtheorem{thm}{Theorem}[section]

 \newtheorem{prop}[thm]{Proposition}
 \theoremstyle{definition}
 \newtheorem{df}[thm]{Definition}
   
 \newtheorem{rem}[thm]{Remark}
 \numberwithin{equation}{section}

\newcommand{\J}{j}
\newcommand{\jj}{j}

\newcommand{\half}{{\textstyle{1\over2}}}

\newcommand{\R}{\mathbb R}%
\newcommand*{\C}{\mathbb C}%
\newcommand*{\Z}{\mathbb Z}%

\newcommand{\gl}{\mathfrak{gl}}
\newcommand{\g}{\mathfrak{g}}
\newcommand{\h}{\mathfrak{h}}

\newcommand{\p}{\mathfrak{p}}

\newcommand{\so}{\mathfrak{so}}

\newcommand{\scal}{{s}}

\makeatletter
\def\operatorname#1{\mathop{\operator@font #1}\nolimits}%
\makeatother
\DeclareMathOperator{\ad}{ad}

\DeclareMathOperator{\End}{End}
\DeclareMathOperator{\Tr}{Tr}
\DeclareMathOperator{\Ad}{Ad}
\DeclareMathOperator{\id}{Id}
\DeclareMathOperator{\Aut}{Aut}

\newcommand{\Span}{\operatorname{Span}}
\newcommand{\Image}{\operatorname{Im}}

\DeclareMathOperator{\Id}{Id}

\newcommand*{\cyclic}{\mathop{\kern0.9ex{{+}\kern-2.2ex\raise-.29ex%
      \hbox{\Large\hbox{$\circlearrowright$}}}}\limits}

\title{Some pseudo-K\"ahler Einstein $4$-symmetric spaces with a ``twin'' special almost complex structure.}

\author{
Michel Cahen$^{(1)}$,  Simone Gutt$^{(1,2)}$,\\  
 Manar Hayyani$^{(3)}$ and Mohammed Raouyane$^{(4)}$\\
\scriptsize{michel.cahen@ulb.be,   
 simone.gutt@ulb.be, manar.hayyani@yahoo.fr, mraouyane@gmail.com}\\
\footnotesize{(1)  D\'{e}partement de Math\'{e}matique, Universit\'{e} Libre de Bruxelles}\\[-7pt]
\footnotesize Campus Plaine, CP 218, Boulevard du Triomphe\\[-7pt]
\footnotesize BE -- 1050 Bruxelles, Belgium\\[-7pt]
\footnotesize Membre de l'Acad\'emie Royale de Belgique.\\
\footnotesize (2) Universit\'e de Lorraine\\ [-7pt]
\footnotesize Institut Elie Cartan de Lorraine, UMR 7502,\\[-7pt]
\footnotesize 3 Rue Augustin Fresnel, 57070 Metz,  France. \\
\footnotesize (3) Universit\'e Moulay Isma\"\i l,\\ [-7pt]
\footnotesize  Facult\'e des Sciences,\\[-7pt]
\footnotesize  Mekn\`es, Maroc.\\
\footnotesize (4) Professeur retrait\'e de l'E.N.S. Universit\'e Mohammed V \\ [-7pt]
\footnotesize  Rabat, Maroc.
}

\date{~}

\begin{document}

\maketitle

\begin{abstract}
On $4$-symmetric symplectic  spaces, invariant almost complex structures -up to sign-  arise in pairs. We exhibit some $4$-symmetric symplectic  spaces, with a pair of  ``natural'' compatible (usually not positive)  invariant almost complex structures,  one of them    being integrable and the other one  being maximally non integrable (i.e. the image of its Nijenhuis tensor at any point is the whole tangent space at that point).  The integrable one defines a pseudo-K\"ahler Einstein metric on the manifold, and the non integrable one  is  Ricci Hermitian (in the sense that the almost complex structure preserves the Ricci tensor of the associated Levi Civita connection) and special in the sense that the associated Chern Ricci form is proportional to the symplectic form. 
 \end{abstract}

\section*{Introduction}

An {\it{almost complex structure}} $\J$  on  a manifold $M$ is  a smooth field of endomorphisms of the tangent bundle whose square is equal to minus the identity. It is said to be {\it{integrable}} if it is induced by a complex structure on $M$; this means that one can locally define complex coordinates on $M$ and that the changes of coordinates are holomorphic; the associated almost complex structure is then given by $\J \frac{\partial}{\partial x^k}= \frac{\partial}{\partial y^k}$ if $x^k$ and $y^k$ are the real and imaginary part of the local complex coordinates $\{z^k=x^k+iy^k \,\vert\, k=1,\ldots,n\}$.\\
The Newlander-Nirenberg theorem asserts that an almost complex structure $\J$ on a manifold $M$  is integrable if and only if its Nijenhuis torsion $N^\J$ vanishes identically.
Recall that the {\it{Nijenhuis   torsion associated to a smooth field $A$ of endomorphisms  of the tangent bundle}} is the tensor of type  
 $(1,2)$ defined by
\begin{equation}
N^A(X,Y):=[AX,AY]-A[AX,Y]-A[X,AY]+A^2[X,Y] \quad \forall X,Y\in \mathfrak{X}(M),
\end{equation}
where $\mathfrak{X}(M)$ is the Lie algebra of ${C}^{\infty}$ vector fields on M.\\
The  {\it{Image distribution}} associated to an almost complex structure $\J$ is 
\begin{equation}\label{def:ImNJ}
(\Image N^\J)_x:=\Span\{ N^\J_x(X,Y)\,\vert\, X,Y \in T_xM \}.
\end{equation}
A way to measure the non integrability of an almost complex structure $\J$ is to look at the dimension of this image distribution; remark that this distribution is smooth but its dimension may vary from point to point.\\
 An almost complex structure $\J$ is called {\it{maximally non integrable}} if $\Image N^\J$ is the whole tangent bundle.\\
 The condition to be maximally non integrable is  open in the space of almost complex structures and generic in high dimension; in a recent preprint \cite{bib:CPS}, R. Coelho, G. Placini, and J. Stelzig show that they exist on any $2n$-dimensional almost complex manifold when $2n\ge 10$.
 We want here to provide such maximally non integrable almost complex structures which are geometrically ``natural".
The first example of a geometrical non integrable almost complex structure  appears in Eells and Salamon  \cite{bib:ES}; it arises on a twistor space when flipping the sign of the vertical part of the standard integrable almost complex structure on this space. We showed in \cite{bib:CGR} that many twistor spaces are  endowed in this way with one integrable and one maximally non integrable almost complex structures.

On a symplectic manifold $(M,\omega)$, 
%(whose dimension will be denoted $2n$ in this paper),
an almost complex structure $\J$  is said to be {\it{ compatible}} with $\omega$ if  the tensor $g_{\J}$ defined by 
\begin{equation}\label{eq:gj}
g_{\J}(X,Y):=\omega(X,\J Y) 
\end{equation}
is symmetric, hence yields an {\it{associated  pseudo-Riemannian metric}}.
A compatible almost complex structure is said to be {\it{positive}} when the associated metric $g_{\J}$ is Riemannian.\\
A triple $(M,\omega, \J)$  with $(M,\omega)$ symplectic and $\J$ a  compatible almost complex structure is equivalent to the data of an {\it{almost pseudo-K\"ahler manifold}} $(M,g,\J)$, i.e. an almost  pseudo-Hermitian manifold 
(which is  a pseudo-Riemannian manifold $(M,g)$ with an almost complex structure $\J$  which is compatible in the sense that the tensor $\omega$ defined by $\omega(X,Y)=g(\J X,Y)$ is skewsymmetric),  with the extra condition  that $d\omega=0$.\\

It is well known that there exist compatible positive almost complex structures on any symplectic manifold, and there have been  various attempts to get procedures to select some  of those.\\
 D. Blair and S. Ianus \cite{bib:BlairIanus} studied  the restriction of the Hilbert functional, usually defined on the space of Riemannian metrics on $M$, to the space of metrics built from positive compatible almost complex structures, i.e.
  $$
\mathcal{F}(\J)=\int_M  \scal^{g_{\J}} \frac{\omega^N}{n!}
$$
where  $\scal^{g_{\J}}$ is the scalar curvature of the Levi Civita connection associated to the metric $g_{\J}$. Extrema of the Hilbert functional are the Einstein metrics.
The extrema of $\mathcal{F}$ are those $\J$'s which are Ricci Hermitian in the following sense :
\begin{df}\label{def:Riccihermitian}
An admissible almost complex structure $\J$ on a symplectic manifold $(M,\omega)$ is said to be {\emph{Ricci Hermitian}} iff
\begin{equation}\label{eq:RicciHerm}
Ric^{g_{\J}}(\J X,\J Y)=Ric^{g_{\J}}(X,Y),
\end{equation}
where $Ric^{g_{\J}}$  is the Ricci tensor of the Levi Civita connection associated to the pseudo-Riemannian metric $g_{\J}$.
\end{df}
This is  automatically satisfied if $\J$ is integrable since the integrability is equivalent to the fact that $\J$ is parallel for the Levi Civita connection.
Examples of non K\"ahlerian triples $(M,\omega,\J)$ satisfying this condition
have been given in 1990 par Davidov et Muskarov \cite{bib:DavidovMusk} on some twistor spaces over Riemannian manifolds.\\

 \noindent V. Apostolov and T. Draghici have introduced in \cite{bib:ApostolovDrag} the notion of  special almost complex structures on a symplectic manifold.
 \begin{df}\label{def:special}
An admissible almost complex structure $\J$ on a symplectic manifold $(M,\omega)$ is said to be {\emph{special}} if its Chern Ricci form is proportional to $\omega$. The Chern Ricci form is the $2$-form on $M$ defined by 
\begin{equation}\label{eq:ChernRicci}
ChernRicci^\J(X,Y):= \Tr \J R^C(X,Y)
\end{equation}
where $R^C$ is the curvature of the Chern connection $\nabla^C$ defined by
\begin{equation}\label{eq:ChernConn}
\nabla^C_XY:= \nabla^{g_{\J}}_XY-\half \J\left(\nabla^{g_{\J}}_X\J\right)Y
\end{equation}
with $\nabla^{g_{\J}}$ denoting the Levi Civita connection associated to the pseudo-Riemannian metric $g_{\J}$. It is the only connection for which $\omega$ and $\J$ are parallel and whose torsion is proportional to $N^\J$.
\end{df}
Positive admissible special almost complex structures were studied by Alberto Della Vedova in \cite{bib:DellaVedova} in a homogeneous context. Homogeneous almost K\"ahler manifolds are also studied by D. V. Alekseevsky and Fabio Podest\`a in \cite{bib:AlekseevskyPodesta}. \\

There exist  two natural almost complex structures $\J^{\pm}$ on the twistor space $\mathcal{T}$ over the hyperbolic space $H_{2n}$. This twistor space can be viewed as an  adjoint orbit in the semisimple Lie algebra $\so(1,2n)$, and thus carries a Kirillov-Kostant-Souriau symplectic structure $\omega$ .
It was shown in \cite{bib:CGGH} that both  $\J^{\pm}$ are compatible with $\omega$, $\J^-$  is positive  and  maximally non integrable whereas $\J^+$ is integrable but  not positive.\\
This twistor space is an example of a $4$-symmetric space and we study in this paper a generalization of this construction on some $4$-symmetric spaces; many examples appeared in  the thesis of Manar Hayyani \cite{bib:theseMH}.
\\

Symplectic and almost complex structures on $k$-symmetric spaces have been  considered, and many results exist on $3$-symmetric spaces. In \cite{bib:Tr}, Maciej Bochenski and  Aleksy Tralle study symplectic structures on $k$-symmetric spaces, which are invariant by the symmetries,  and a list of all symplectic 3-symmetric manifolds with simple groups of transvections is given. This extends our previous work about symmetric symplectic spaces \cite{bib:BCG}  and the results obtained by Pierre Bieliavsky \cite{bib:Biel}. \\
Cecilia Ferreira studies in \cite{bib:ferreira}  some necessary and some sufficient conditions for the integrability of a canonical  almost complex structure on a k-symmetric space. These always exist when $k$ is odd.
J.A. Jim\'enez \cite{bib:Jimenez} has given a classification of compact simply connected Riemannian $4$-symmetric spaces; these spaces are homogeneous for a  connected compact semisimple Lie group with an automorphism  of order four.  Geometrically, they can be regarded as fiber bundles over Riemannian $2$-symmetric spaces with totally geodesic fibers isometric to a Riemannian $2$-symmetric space. \\
Vitaly Balashchenko studies in \cite{bib:Balash} canonical distributions on Riemannian homogeneous $k$-symmetric spaces; those are distributions associated to canonical affinor structures (an affinor structure is a field of endomorphisms $A$ which is  a product structure, i.e. $A^2=\id$,
 an almost complex structure, i.e. $A^2=-\Id$, a $f$-structure, i.e. $A^3-A=0$, or a $h$-structure, i.e. $A^3+A=0$).\\

Here we study a class of $4$-symmetric spaces which carry two associated invariant almost complex structures, which are natural in a sense that we define. One of them is integrable and we state conditions for the other one to be maximally non integrable.\\
Some of  those spaces carry a symplectic structure for which 
the almost complex structures are compatible, in general not positive. 

We get families of examples of Einstein pseudo-K\"ahler manifolds, admitting another natural almost complex structure which is maximally non integrable, Ricci Hermitian in the sense of definition \ref{def:Riccihermitian}, and special in the sense of definition \ref{def:special}.

In section \ref{section:sym}, we recall basic facts about  $k$-symmetric spaces and invariant structures. In section \ref{section:constr}, we construct a class of $4$-symmetric spaces with two ``twin''  invariant almost complex structures $\J^\pm$. We indicate  conditions under which  a symplectic  $2$-form $\omega$  can be defined, for which $\J^\pm$ are compatible. 
The examples  given in section \ref{section:examples} are summarized in :
\begin{prop}
Let  $G$ be a connected Lie subgroup of $Gl(m,\R)$, with $m=k+2n$, whose Lie algebra $\g$ is stable under $\ad \rho$ where $\rho=\left( \begin{smallmatrix}
0 & 0\cr
0 & J_{2n}
\end{smallmatrix}\right)$ with $J_{2n}=\left(\begin{smallmatrix}
0 & -\Id_n\cr
\Id_n & 0
\end{smallmatrix}\right)$. Denote by $\sigma$ the automorphism of $\g$ of order $4$ given by
$\sigma=\exp (\half \pi \ad \rho)$.\\
$G$ is 
stable by the automorphism
$
\tilde\sigma:G\rightarrow G :
g\mapsto RgR^{-1}$ for $R:= \left( \begin{smallmatrix}
\Id_k & 0\cr
0 & J_{2n}
\end{smallmatrix}\right)$ and the homogeneous space $G/G^{\tilde\sigma}_0$, where $G^{\tilde\sigma}_0$ is the connected component of the subgroup of elements fixed by $\tilde\sigma$, has a natural structure of $4$-symmetric space. \\
The tangent space at the base point identifies with $\g^\sigma_{-1}\oplus \g^{\sigma^2}_{-1}$ where $\g^\tau_{\lambda}$ denotes  the $\tau$ eigenspace of eigenvalue $\lambda$ in $\g$.
Two natural  invariant almost complex structures $\J^\pm $ on $G/G^{\tilde\sigma}_0$ are defined by their value at the base point, $J^\pm\in \End(\g^\sigma_{-1}\oplus \g^{\sigma^2}_{-1})$, given by
$$
J^\pm\vert_{\g^{\sigma^2}_{-1}}:=\sigma\vert_{\g^{\sigma^2}_{-1}}\qquad J^\pm\vert_{\g^\sigma_{-1}}=\pm \exp (\frac14 \pi \ad \rho);$$
$\J^+$ is always integrable; the image of the Nijenhuis torsion of $\J^-$ is the $G$-invariant distribution whose value at the base point is given by 
$[\g^{\sigma^2}_{-1},\g^\sigma_{-1}]+[\g^{\sigma^2}_{-1},\g^{\sigma^2}_{-1}]\cap \g^\sigma_{-1}.$\\
The $G$-invariant $2$-form $\omega$ on $G/G^{\tilde\sigma}_0$ whose value at the base point is given by 
$$
\tilde{\Omega}(X,Y):= \Tr(\rho[X,Y]) \qquad \forall X,Y \in \g^\sigma_{-1}+\g^{\sigma^2}_{-1}
$$
is an invariant symplectic structure, if and only if it is non degenerate, i.e. iff $
g' :\g^{\sigma^2}_{-1}\times \g^{\sigma^2}_{-1}\rightarrow \R :X,Y \mapsto \Tr(XY)$ and $ \beta^v : \g^\sigma_{-1} \times \g^\sigma_{-1} \rightarrow \R : X,Y \mapsto \Tr(XY)$ are non degenerate; this is automatically true if $\g$ is simple. The almost complex structures $\J^\pm$ are compatible with $\omega$, in general not positive. When $\omega$ is non degenerate, one considers the associated metrics, $g^\pm(X,Y)=\omega(X,\J^\pm Y)$. The Ricci tensors of  the associated Levi Civita connections are $\J^\pm$ hermitian and, under some conditions on $\g$, the structure $\J^-$ is special. In particular, we have the following examples.\\

The space $Sl(k+2n,\R)/S(Gl(k,\R)\times Gl(n,\mathbb{C}))$, where $Gl(n,\mathbb{C})$ is the subset of elements in $Gl(2n,\mathbb{R})$ commuting with $J_{2n}$, 
with the pair $(\omega,\J^+)$, is a  pseudo-K\"ahler $4$-symmetric space which is Einstein
$(Ric^{g^+}=(k+n)\, g^+)$;
with the pair $(\omega,\J^-)$, it is an  almost-pseudo-K\"ahler $4$-symmetric space,
$\J^-$ is maximally non integrable and   special $(ChernRicci^{j^-}=2(n-k)\, \omega)$.\\

The space $SO(k+2n,\R)/(SO(k,\R)\times U(n))$, where $U(n)$ is the subset of elements in $SO(2n,\R)$ commuting with $J_{2n}$,  with the pair $(-\omega,\J^+)$, is a  K\"ahler $4$-symmetric space which is Einstein $(Ric^{g^+}=\half(k+n-1)g^+$, $g^+$ negative definite$)$ ; with the pair $(\omega,\J^-)$, it is 
an almost pseudo-K\"ahler $4$-symmetric space,  $\J^-$ is maximally non integrable and special $(ChernRicci^{j^-}=(n-1-k)\omega)$.
The case where $k=1$ corresponds to the twistor space $SO(1+2n,\R)/ U(n)$ on the sphere $SO(1+2n,\R)/ SO(2n,\R)$.\\

The space $ SO_0(k,2n,\R)/(SO(k,\R)\times U(n))$, with the pair $(\omega,\J^+)$, is a pseudo-K\"ahler $4$-symmetric space which is Einstein $(Ric^{g^+}=\half(k+n-1)g^+)$; with the pair $(\omega,\J^-)$, it is 
an almost K\"ahler $4$-symmetric space,  $\J^-$ is maximally non integrable and special $(ChernRicci^{j^-}=(n-1-k)\omega)$.
The case where $k=1$ corresponds to the twistor space $SO_o(1,2n,\R)/ U(n)$ on the hyperbolic space  $SO_o(1,2n,\R)/SO(2n,\R)$. \\

The  space $U(k'+n)/\left(U(k')\times U(n)\right)$ with the pair $(-\omega, \J^+)$ is a K\"ahler symmetric space which is Einstein $(Ric^{g^+}=\half(k'+n)g^+$,  $g^+$ negative definite$)$.
The space  $ U(k',n)/\left(U(k')\times U(n)\right)$ with the pair $(\omega, \J^+)$ is a  K\"ahler symmetric space which is Einstein $(Ric^{g^+}=\half(k'+n)g^+)$. 
Both  spaces carry only one natural almost complex structure (up to sign), $\J^+$, defined by $\sigma$ on $\p=\g^{\sigma^2}_{-1}$. \\

The space $Sp(\R^{2(k'+n)},\widetilde\Omega)/(Sp(\R^{2k'},\Omega_k)\times U(n))$, 
{{where  $\widetilde\Omega= \left(\begin{smallmatrix}
\Omega_{2k'}&0\\0&\Omega_{2n}
\end{smallmatrix}\right)$ with $\Omega_{2r}:=-J_{2r}$ and  
where $U(n)$ is the subset of elements in $Sp(\R^{2n},\Omega_{2n})$ commuting with $J_{2n}$,}}
with  $(\omega,\J^+)$, is a pseudo-K\"ahler $4$-symmetric space which is Einstein
$(Ric^{g^+}=\half(2k'+n+1)g^+)$;  with the pair $(\omega,\J^-)$, it is an almost pseudo-K\"ahler $4$-symmetric space, $\J^-$ is maximally non integrable and  special $( ChernRicci^{j^-}=(n+1-k)\omega)$.\\

For all those spaces, the symplectic structure coincides with the Kirillov-Kostant-Souriau symplectic form on the coadjoint orbit of the element $\rho^\flat: \g \rightarrow \R : X \mapsto \Tr \rho X$ in $\g^*$. 
\end{prop}

 \section*{Acknowledgement} This work benefited from the ``Excellence of Science (EoS)" grant number 30950721 ``Symplectic techniques in differential geometry", funded by  the FWO/F.R.S-FRNS.

\section{$k$-symmetric spaces}\label{section:sym}

 The notion of $k$-symmetric space appeared in the late sixties in the works of Gray, Wolf, Ledger, and  Obata \cite{Gray, WG, Ledger, L-O}, as a generalization of the notion of symmetric space which corresponds to $k=2$.
General results are summarized in a monograph written by Kowalski \cite{Kow} in 1980. 
\begin{df}A {\emph{$k$-symmetric space}} (with $k$ an integer $\ge 2$) is a pair $(M,S)$, where $M$ is a smooth manifold and 
$S:M\times M \rightarrow M : (x,y)\mapsto S(x,y)=:s_xy$
is a smooth map such that \\
- for each $x\in M$, $s_x$ -which will be called  the symmetry at $x$- is a  diffeomorphism of order $k$ of $M$  (i.e. $s_x^k=\Id$ and $k$ is the smallest positive integer with that property),
 which fixes $x$,  and the differential of  $s_x$ at the point  $x$ has no non-zero fixed vector (which implies that $x$ is an isolated fixed point); \\
- $s_x\circ s_y =s_{s_x y}\circ s_x $ for all $ x,y \in M$.
\end{df}
An {\emph{automorphism of a $k$-symmetric space}} $(M,S)$ is a diffeomorphism $\varphi:M\rightarrow M$ such that $
\varphi\left(s_x(y)\right)=s_{\varphi(x)}\varphi(y)\quad \forall x,y \in M.
$ Each symmetry is an automorphism.\\
It is known that, if $M$ is connected, the group of automorphisms $\Aut(M,S)$ is a Lie group which acts transitively on $M$; we shall only consider $k$-symmetric spaces for which  $\Aut(M,S)$ acts transitively.  Choosing a point $p_0\in M$, one defines the automorphism ${\tilde\sigma}$ of $\Aut(M,S)$ given by 
$g\mapsto {\tilde\sigma}(g):=s_{p_0}\circ g \circ (s_{p_0})^{-1}$.
The stabilizer $H$ of $p_0$ in $\Aut(M,S)$ is contained in the subgroup $\Aut(M,S)^{\tilde\sigma}$ of elements fixed by ${\tilde\sigma}$, and the connected component of $\Aut(M,S)^{\tilde\sigma}$ is contained in $H$.
\begin{df}
A {\emph{$k$-symmetric  triple}} (also called \emph{ a homogeneous $k$-symmetric space} for the Lie group $G$) is a triple $(G,{\tilde\sigma},H)$ where $G$ is a Lie group, ${\tilde\sigma}$ is an automorphism of $G$ of order $k$ and $H$ is a subgroup of $G$ such that 
\begin{equation}\label{eq:H}
G^{\tilde\sigma}_\circ \subset H \subset G^{\tilde\sigma} 
\end{equation}
where  $G^{\tilde\sigma}$ is the subgroup of elements of $G$ which are fixed by ${\tilde\sigma}$ and $G^{\tilde\sigma}_\circ$ its connected component of the identity.
\end{df}
A $k$-symmetric triple  defines a $k$-symmetric space $(M,S)$ where $M$ is the homogeneous space $M=G /H$,
 and  $S:M\times M\rightarrow M$   is defined by      
$$
 s_{\pi(g)}\pi(g'):=\pi(g {\tilde\sigma} (g^{-1}g')),\quad \forall g,g' \in G,
$$  
where $\pi:G\rightarrow G/H$ denotes the canonical projection. Remark that $G$ acts by automorphisms of this symmetric space.\\
Hence, a  $k$-symmetric triple defines a $k$-symmetric space in a unique way, but a connected $k$-symmetric space can be associated to different $k$-symmetric triples.
Since we study invariant structures, we only deal with homogeneous spaces and the precision of the group of invariance considered will be of crucial importance;  we shall thus use the description of a (homogenous) $k$-symmetric space in terms of a chosen $k$-symmetric triple.

A $k$-symmetric  triple defines  a pair $(\g,\sigma)$ where $\g$ is the Lie algebra of $G$ and where 
$\sigma:\g \rightarrow \g$ is the automorphism of $\g$ of order $k$ defined by the differential of ${\tilde\sigma}$ at the neutral element $e\in G$. The Lie algebra $\h$ of the group $H$ is always equal to the subalgebra $\g^\sigma$ of elements which are fixed by $\sigma$.

\begin{df}\label{def:algksym} A {\it{$k$-symmetric algebra}} (where $k$ is an integer $\ge 2$) is a pair $(\g,\sigma)$, where $\g$ is a Lie algebra and $\sigma$ is an automorphism of $\g$ of order $k$.
 \end{df}
 
A $k$-symmetric algebra defines a $k$-symmetric triple $( G,{\tilde\sigma},H)$ for any   Lie group $G$ with Lie algebra  $\g$  which admits an automorphism $\tilde\sigma$ lifting  $\sigma$ (for instance the simply connected group) and  for any  $H$ satisfying relation \ref{eq:H} (for instance $H= G^{\tilde\sigma}_\circ$).\\
Observe that, for any k-symmetric space corresponding to a $k$-symmetric algebra $(\g,\sigma)$, the tangent space to $M=G/H$ at the base point $p_0=\pi(e)$ identifies to $\g/\g^\sigma$ and the differential of the symmetry $s_{p_0}$ at the base point identifies with the linear map on $\g/\g^\sigma$ induced by $\sigma$. 

\begin{df}
 An {\it{invariant almost complex structure}} (resp. invariant symplectic structure, invariant pseudo-riemannian structure,..) on a homogeneous $k$-symmetric space described by a {$k$-symmetric  triple} $(G,\tilde\sigma, H)$ is an almost complex structure (resp. symplectic structure, pseudo-riemannian structure,..) on $M=G/H$ which is invariant by the symmetries and by the action of $G$. \\
Such a structure is completely determined by its value at the base point, and that restriction must be invariant under the action of $H$ and under the symmetry at $p_0$. Observe that the differential at the base point of the action of an element $h\in H$
is the action, still denoted ${\Ad h}$, induced by $\Ad h$ on the quotient $\g/\g^\sigma$.
The differential at the base point of the symmetry at the base point is given by the action,  still denoted ${\sigma}$, induced by $\sigma$ on $\g/\g^\sigma$.
The action of $\sigma$ and the action of $\Ad h$ commute since $H\subset G^{\tilde\sigma}$.
If $H$ is connected, a structure is invariant iff its value at the base point is invariant under  the action induced on the quotient $\g/\g^\sigma$ by $\sigma$ and by $\ad X$ for all $X\in \g^\sigma$.\\
An {{invariant almost complex structure}} $\J$  corresponds thus bijectively to a 
$$J:\g/\g^\sigma\rightarrow  \g/\g^\sigma$$ such that $J\circ \sigma=\sigma \circ J$ and $J\circ {\Ad h}={\Ad h} \circ J$ for all $h\in H$.
\end{df}

To study the integrability of an invariant almost complex structure on a homogeneous $k$-symmetric space, we shall use the following description:
\begin{prop}(\cite{KN}, thm 6.4, page 217)\label{prop:Nijinv}
Let $\J$ be a $G$-invariant almost complex structure on a homogeneous space $M=G/H$, let $J:\g/\h\rightarrow  \g/\h$ be its value at the base point $p_0=\pi(e)$, with $\pi: G\rightarrow G/H$ the canonical projection. Let $\hat J:\g\rightarrow \g$ be any linear map such that
$\pi_*\hat J(X)=J(\pi_*X)$ with $\pi_*$ the differential of $\pi$ at the neutral element $e$, i.e.
the canonical projection $\pi_*:\g\rightarrow \g/\h$. Then $$N^\J_{p_0}(\pi_*X,\pi_*Y)=\pi_*\tilde N^{\hat J}(X,Y)$$
where $\tilde N^{\hat J}(X,Y):=[\hat J X,\hat J Y]-\hat J[\hat J X,Y]-\hat J[X,\hat J Y]+\hat J^2 [X,Y]$ for all $X,Y \in \g$.
\end{prop}

We are particularly interested to build  homogeneous $4$-symmetric (almost)-pseudo-K\"ahler manifolds, hence study invariant complex structures $\J$ on a $4$-symmetric 
symplectic space   which are compatible with the invariant symplectic $2$-form $\omega$.

 \begin{df}\label{def:sympl} A {\it{homogeneous  $k$-symmetric symplectic space}} is a $G$-homogeneous $k$-symmetric space 
 endowed with a symplectic $2$-form  which is $G$-invariant and 
 invariant by the symmetries; it is defined by a {\it{$k$-symmetric symplectic quadruple}} $(G,\tilde\sigma,H,\tilde\Omega)$ with $(G,\tilde\sigma,H)$ a $k$-symmetric triple and   
 $$\tilde\Omega:\g/\g^\sigma\times \g/\g^\sigma \rightarrow \R$$
a non degenerate bilinear skewsymmetric map such that
 \begin{itemize}
 \item $\tilde\Omega (\pi_*(\Ad h X), \pi_*(\Ad h Y))=\tilde\Omega(\pi_*(X),\pi_*(Y))$ for all $X,Y \in \g$, for all $h\in H$;
 \item $\tilde\Omega (\pi_*(\sigma  X), \pi_*(\sigma Y))=\tilde\Omega(\pi_*(X),\pi_*(Y))$ for all $X,Y \in \g$;
 \item $\cyclic_{XYZ}\tilde\Omega (\pi_*([X,Y]), \pi_*(Z))=0$ for all $X,Y,Z \in \g$,
 \end{itemize}
$\cyclic_{XYZ}$ denoting the sum over cyclic permutations and  $\pi_*: \g\rightarrow \g/\g^\sigma$ being, as before,  the canonical projection.
 \end{df}

Gray and Wolf have introduced in \cite{WG} the notion of a canonical field of endomorphisms:
 \begin{df}
 A field of endomorphisms  on a $k$-symmetric space described by a {$k$-symmetric  triple} $(G,\tilde\sigma, H)$ is said to be {\it{canonical}} if it is invariant and if its value at the base point
is given by the action on $\g/\g^\sigma$ of a polynomial in $\sigma$.\\
\noindent Remark that any polynomial in $\sigma$ is automatically invariant under  the action induced
on the quotient $\g/\g^\sigma$ by $\sigma$ and by $\Ad h$ for all $h\in H$ so yields a canonical field of endomorphisms.
\end{df}

Any endomorphism of order $k$ of a finite dimensional vector space  being semisimple,  one can identify   the complexification of the tangent space at the base point, $\left(\g/\g^\sigma\right)^\mathbb{C}$, to the sum of all eigenspaces for $\sigma$, corresponding to eigenvalues (which are of the form $\lambda=e^{2i\pi r/k}$  since $\sigma^k=\Id$) which are different from $1$.

To have a canonical almost complex structure on a $k$-symmetric space,  $-1$ can not be   an eigenvalue for $\sigma$. This is of course automatic if $k$ is odd. In this paper we want to study  almost complex structures on $4$-symmetric spaces and we have to go beyond this canonical condition. 

\begin{df}\label{def:naturalJ}
Let $(G,\tilde\sigma,H)$ be a $4$-symmetric triple and write the decomposition of the Lie algebra $\g$ :
$$\g=\g^{\sigma^2}+\g^{\sigma^2}_{-1}=\g^\sigma+\g^\sigma_{-1}+\g^{\sigma^2}_{-1},$$ 
where $\g^{\sigma^k}_{\lambda}$ is the eigenspace for $\sigma^k$ of eigenvalue $\lambda$. One  identifies, via $\pi_*$, the tangent space to $G/H$ at the base point $eH$ with $\g^\sigma_{-1}+\g^{\sigma^2}_{-1}$. Any invariant almost complex structure $j$ on $G/H$ is defined by 
$$
J:\g^\sigma_{-1}+\g^{\sigma^2}_{-1} \rightarrow \g^\sigma_{-1}+\g^{\sigma^2}_{-1}
$$
such that $J_{\vert{\g^{\sigma^2}_{-1}}}=\hat\sigma: \g^{\sigma^2}_{-1}\rightarrow \g^{\sigma^2}_{-1}$ and $J_{\vert_{\g^\sigma_{-1}}}=\tau:\g^\sigma_{-1}\rightarrow \g^\sigma_{-1}$ satisfy
$$
(\hat\sigma)^2=\sigma^2_{\vert{\g^{\sigma^2}_{-1}}}=-\Id_{\vert{\g^{\sigma^2}_{-1}}} , \quad \hat\sigma\circ \sigma_{\vert{\g^{\sigma^2}_{-1}}}= \sigma_{\vert{\g^{\sigma^2}_{-1}}} \circ \hat\sigma, \quad \hat\sigma\circ \Ad(h)= Ad(h)\circ \hat\sigma \quad \forall h\in H;
$$
$$
\tau^2=\sigma_{{\g^\sigma_{-1}}} =-\Id_{{\g^\sigma_{-1}}}   \quad \textrm{and }\quad \tau\circ \Ad(h)= Ad(h)\circ \tau \quad \forall h\in H.
$$
An invariant almost complex structure on a $4$-symmetric space is said to be {\emph{natural}} if $\hat\sigma=\sigma$, i.e. if $J_{\vert{\g^{\sigma^2}_{-1}}}=\sigma_{\vert \g^{\sigma^2}_{-1}}$.\\
When $\g^\sigma_{-1}\neq 0$, invariant almost complex structures on $4$-symmetric spaces  arise  in pairs (called {\emph{twins}}) $j^{\pm}$ corresponding to $J^\pm $ defined by
$$
J^\pm\vert_{\g^{\sigma^2}_{-1}}=\hat\sigma \quad \textrm{and }\quad J^\pm\vert_{\g^\sigma_{-1}}=\pm\tau.
$$

\end{df}

\begin{rem}\label{rem:transvections}  
Any homogeneous $4$-symmetric space $(G,\tilde\sigma,H)$ fibers over the symmetric space $G/G^{{\tilde\sigma}^2}$ since $(\tilde\sigma^2)^2=\id$ and $H\subset G^{{\tilde\sigma}}\subset G^{{\tilde\sigma}^2}$.\\ Observe that $\g^{{\sigma}^2}=\g^\sigma + \g^\sigma_{-1}$ and that $G$ is contained in the group of automorphisms of the symmetric space $G/G^{{\tilde\sigma}^2}$. The Lie algebra $\g'$ of the transvection group for this symmetric space, which is given by 
$$
\g'=\p +[\p,\p] \quad \textrm{ with } \p=\g^{\sigma^2}_{-1},
$$
is an ideal in $\g$.\\
The fiber of the projection $G/H\rightarrow G/G^{{\tilde\sigma}^2}$ identifies with the symmetric space $G^{{\tilde\sigma}^2}/H$ defined by the triple $(G^{{\tilde\sigma}^2},\tilde\sigma,H)$. \\ 

A natural almost complex structure $\J$ on $G/H$ induces an invariant almost complex structure on the symmetric space $G/G^{{\tilde\sigma}^2}$ if and only if $\sigma$ commutes with $\Ad k$ for any $k\in G^{{\tilde\sigma}^2}$; in that case
the fibration $G/H \rightarrow G/G^{{\tilde\sigma}^2}$ is a pseudo-holomorphic map.\\
The fiber  $G^{{\tilde\sigma}^2}/H$ always carries an invariant complex structure $\J^v$ induced by $\tau$
since the tangent space to the fiber at $eH$ identifies to $\g^{{\sigma}^2}/\h=\g^{{\sigma}^2}/\g^{{\sigma}}\simeq\g^\sigma_{-1}$; it is integrable since $[\g^\sigma_{-1},\g^\sigma_{-1}]\subset \g^\sigma$.
\end{rem}

\begin{df}
A {\emph{homogeneous $4$-symmetric natural (almost)-pseudo-K\"ahler manifold}} is
a homogeneous $4$-symmetric space $(G,\tilde\sigma,H)$, endowed with a natural almost complex structure $\J$, defined by $\tau$ as in definition \ref{def:naturalJ},  and  with a symplectic structure $\omega$ on $G/H$ invariant by $G$ and by the symmetries, defined by a $\tilde\Omega$ as in definition \ref{def:sympl}, such that $\J$ is compatible with it. 
Since the tangent space to $G/H$ at the base point identifies via $\pi_*$ with $\g^\sigma_{-1}+\g^{\sigma^2}_{-1}$, the data of $\tilde\Omega$ with the compatibility with $J$ implies the compatibility with $J^\pm$ and is equivalent to the data of two non degenerate skewsymmetric bilinear maps
$$
\tilde\Omega^v:\g^\sigma_{-1}\times \g^\sigma_{-1}\rightarrow \R \qquad \qquad 
\tilde\Omega':\g^{\sigma^2}_{-1}\times \g^{\sigma^2}_{-1}\rightarrow \R
$$
such that 
\begin{eqnarray}
&&\tilde\Omega^v(\tau X,\tau X')=\tilde\Omega^v(X,X')=\tilde\Omega^v(\Ad h X, \Ad h X')\\ 
&&\tilde\Omega'(\sigma Y,\sigma Y')=\tilde\Omega'(Y,Y')=\tilde\Omega'(\Ad h Y, \Ad h Y')\\
&&\tilde\Omega' ([X,Y], Y')+\tilde\Omega' (Y,[X,Y'])+\tilde\Omega^v (\half\left([Y,Y']-\sigma([Y,Y'])\right), X)=0 \label{eq=omegav}
\end{eqnarray}
for all  $X,X' \in \g^\sigma_{-1} , \quad  h\in H$ and $Y,Y' \in \g^{\sigma^2}_{-1}$.
\end{df}
The fiber $G^{{\tilde\sigma}^2}/H$ is then a symmetric   pseudo-K\"ahler  manifold with the symplectic structure defined by $\tilde\Omega^v$ and the complex structure defined by $\tau$.\\
The basis of the fibration, i.e. the symmetric space $G/G^{{\tilde\sigma}^2}$, is endowed with an invariant symplectic structure defined by $\tilde\Omega'$ iff $\tilde\Omega'( X, X')=\tilde\Omega^v(\Ad k X, \Ad k X')$ for all $k\in G^{{\tilde\sigma}^2}$.  This condition, in view of equation \ref{eq=omegav} and the fact that $\g^\sigma_{-1}\subset \g^{\sigma^2}$, implies that $[\g^{\sigma^2}_{-1},\g^{\sigma^2}_{-1}]\subset \g^{\sigma}$.

 \section{A class of $4$-symmetric spaces}\label{section:constr}

Let $D$ be a   derivation of a Lie algebra $\g$ such that $\exp 2 \pi D=\Id$. Remark that this implies that $D$ is semisimple and that all eigenvalues of $D$ are contained in $i\Z$.
Assume also that $\exp t D\neq\Id$ for all $0<t<2\pi$. Let
$$\sigma=\exp \frac12 \pi D.$$
Clearly $\sigma$ is an automorphism of order $4$ and we consider the $4$-symmetric algebra $(\g,\sigma)$. We have, as above,  the splitting
$$
\g=\g^\sigma+\g^\sigma_{-1}+\p \quad\textrm{ with } \quad  \p=\g^{\sigma^2}_{-1}\quad \textrm{ i.e. } \quad\sigma^2\vert_{\p}=-\Id\vert_{\p}.
$$
 We identify $\g/\g^\sigma$ to $\g^\sigma_{-1}+\p$ and we define, as in definition \ref{def:naturalJ}, two natural  complex structures $J^\pm$ on $\g^\sigma_{-1}+\p$ by
$$
J^\pm\vert_{\p}:=\sigma\vert_{\p}\qquad J^\pm\vert_{\g^\sigma_{-1}}=\pm \exp \frac14 \pi D.
$$
They are different iff $\g^\sigma_{-1}\neq \{ 0\}$.
For these to define invariant almost complex structures on $G/G^{\tilde\sigma}_0$, they have to commute with 
$\sigma$ -which is obvious- and  to commute with the action of $\ad X$ for each $X \in  \g^\sigma$;  so  one has only to check if $\exp \frac14 \pi D$ commutes with $\ad X$
on $\g^\sigma_{-1}$ for all $X\in \g^\sigma$. Clearly 
$$
\g^\sigma=\g^\sigma_{1}=\left(\bigoplus_{{m\in\Z}}\, \g^D_{4mi}\right)\cap \g \qquad \g^\sigma_{-1}=
\left(\bigoplus_{{m\in\Z}}\, \g^D_{(4m+2)i}\right)\cap \g
$$
where $\g^D_{mi}$ is the eigenspace of $D$ of eigenvalue $mi$ in the complexified Lie algebra $\g^{\mathbb{C}}$. We have, for all $X\in \g^D_{4mi}$ and $Y\in \g^D_{(4m'+2)i}$
\begin{eqnarray*}
(\ad X\circ \exp \frac14 \pi D )\,Y &=& \left(\exp(m'+\frac12)i\pi\right)\,  [X,Y]=(-1)^{m'}i[X,Y]\\
(\exp \frac14 \pi D\circ \ad X) Y&=& \left(\exp(m'+m+\frac12)i\pi\right) \,[X,Y]=(-1)^{m'+m}i[X,Y]  .
\end{eqnarray*}
The structures $J^\pm$ define thus invariant almost complex structure on $G/G^{\tilde\sigma}_0$ if $4mi$ is not an eigenvalue of $D$ for $m$ odd.\\

We  consider the particular case where $D$ is a semisimple derivation of $\g$ whose eigenvalues are precisely $0, i,-i, 2i$ and $-2i$, and, as above,  $\sigma=\exp \frac12 \pi D$.
In that case 
$$\g^\sigma=\g\cap\g^D_0,\quad  \g^\sigma_{-1}=\g\cap\left(\g^D_{2i}\oplus \g^D_{-2i}\right) \, \textrm{ and } \,\p=\g\cap(\g^D_{i}\oplus \g^D_{-i}).$$
The invariant almost complex structures on $G/G^{\tilde\sigma}_0$ defined by $
J^\pm\vert_{\p}:=\sigma\vert_{\p}$ and $J^\pm\vert_{\g^\sigma_{-1}}=\pm \exp \frac14 \pi D
$ have  Nijenhuis torsions which are $G$-invariant tensors. 
The maps $J^\pm$ extend $\C$-linearly to $T_{p_0}(G/G^{\tilde\sigma}_0)^\C=\g^D_{i}\oplus \g^D_{-i}\oplus \g^D_{2i}\oplus \g^D_{-2i}$ as 
$$
J^\pm\vert\g^D_{i}=i\Id\vert\g^D_{i}\quad J^\pm\vert\g^D_{-i}=-i\Id\vert\g^D_{-i}\quad J^\pm\vert\g^D_{2i}=\pm i\Id\vert\g^D_{2i}\quad J^\pm\vert\g^D_{-2i}=-\pm i\Id\vert\g^D_{-2i}.
$$
Using proposition \ref{prop:Nijinv}, and extending ${N}^{\J^\pm}_{p_0}$ to the complexified tangent space at the base point, we have 
$$
{N}^{\J^\pm}_{p_0}(\pi_*X,\pi_* Y)=\left\{\begin{array}{ll} 0 & \textrm{ for } X,Y \in \g^D_{2i}\, \textrm{ or }\, X,Y\in \g^D_{-2i}\\[1mm]
            0 & \textrm{ for } X \in \g^D_{2i},Y \in \g^D_{-2i} \\[1mm] 
            0   &  \textrm{ for } X \in \g^D_{i},Y \in \g^D_{-i}   \\[1mm] 
             0   &  \textrm{ for } X \in \g^D_{i},Y \in \g^D_{2i}\, \textrm{ or }\, X \in \g^D_{-i},Y \in \g^D_{-2i}  \\[1mm] 
            (-2 \pm 2)\pi_*[X,Y] & \textrm{ for } X,Y \in \g^D_{i} \, \textrm{ or }\, X,Y\in \g^D_{-i}\\[1mm]
            (\pm 2-2)\pi_*[X,Y] & \textrm{ for } X \in \g^D_{i},Y \in \g^D_{-2i}\, \textrm{ or }\,X \in \g^D_{2i},Y \in \g^D_{-i}\\[1mm]
              \end{array}\right.
$$
so that  ${N}^{\J^+}_{p_0}(\pi_*X,\pi_* Y)=0,\,\, {N}^{\J^-}_{p_0}(\pi_*X,\pi_* Y)=-4\pi_*[X,Y]$ for all $X,Y \in \g^\sigma_{-1}\oplus \p$.
Hence we have 
\begin{prop}\label{prop:D}
Let $D$ be a semisimple  derivation of a Lie algebra $\g$ whose eigenvalues are precisely $0, i,-i, 2i$ and $-2i$. Let $\sigma=\exp \frac12 \pi D$. Let $G$ be a Lie group with Lie algebra $\g$ such that $\sigma$ lifts to an automorphism $\tilde\sigma$  of $G$. Let $\J^\pm$ be the two natural  almost complex structure  on the $4$-symmetric space $G/G^{\tilde\sigma}_0$ defined, identifying the tangent space at the base point  to $\g^\sigma_{-1}+\p$, with $\p=\g^{\sigma^2}_{-1}$, by
$$ 
J^\pm\vert_{\p}:=\sigma\vert_{\p}\qquad J^\pm\vert_{\g^\sigma_{-1}}=\pm \exp \frac14 \pi D.
$$
Then $\J^+$ is always integrable. The image of the Nijenhuis torsion of $\J^-$ is the $G$-invariant distribution whose value at the base point is given by 
$$
\pi_*\left([\g^\sigma_{-1}+\p,\g^\sigma_{-1}+\p]\right)=\pi_*\left([\p,\g^\sigma_{-1}]+[\p,\p]\cap \g^\sigma_{-1}\right).
$$
Thus $\J^-$ is maximally non integrable iff $[\p,\g^\sigma_{-1}]= \p$ and $[\p,\p]\cap \g^\sigma_{-1}=\g^\sigma_{-1}$.
\end{prop}
\begin{rem}\label{rmkintegr1}
The equality $[\p,\p]\cap \g^\sigma_{-1}=\g^\sigma_{-1}$  will be true iff $\p +[\p,\p]+\g^\sigma=\g$,
in particular (since  $\p +[\p,\p]$ is an ideal in $\g$ by remark \ref{rem:transvections}) if $\g$ is  simple, or if $\g$ is reductive with $\p$ intersecting each simple factor  and $\sigma=\Id$ on the center.
\end{rem}
%%%%%%%%%%%%%%%%%%%%%%%%%%%%%%%%%
%%%%%%%%%%%%%%%%%%%%%%%%%%%%%%%%%%
%%%%%%%%%%%%%%%%%%%%%%%

\subsection{Construction}\label{section:ex}

\subsubsection{The derivation, the automorphism, and the almost complex structures}

Let  $\rho$ be the $(k+2n)\times (k+2n)$ matrix, with $k$ and $n$ positive integers, defined by
\begin{equation}\label{derfR}
\rho= \begin{pmatrix}
0 & 0\cr
0 & J_{2n}
\end{pmatrix} \qquad \textrm{ where } \qquad J_{2n}=\begin{pmatrix}
0 & -\Id_n\cr
\Id_n & 0
\end{pmatrix}
\end{equation}
and let $D$ be the derivation of $\gl(m,\R)$ ,  with $m=k+2n$, defined by
$$
D=\ad \rho \qquad \textrm{ so that }\qquad  D\begin{pmatrix}
A & B\cr
B' & C
\end{pmatrix} = \begin{pmatrix}
0 & -BJ_{2n}\cr
J_{2n}B' & [J_{2n},C]
\end{pmatrix}.
$$
Clearly $\rho$ is semisimple with eigenvalues $0,i$ and $-i$, so that $D$ is a semisimple derivation with eigenvalues $0,i,-i,2i$ and $-2i$ and we can apply the results above
for any subalgebra $\g$ of $\gl(m,\R)$ which is stable under $D$,  and so that $\g^\C$   intersects 
the $2i$ eigenspace of $D$. Observe that
\begin{equation}
\sigma= \exp \frac12 \pi D= \exp \frac12 \pi \ad \rho = \Ad \exp (\frac12 \pi \rho) =\Ad \begin{pmatrix}
\Id_k & 0\cr
0 & J_{2n}
\end{pmatrix},
 \end{equation}
and we shall denote by $R$ the matrix $R=\begin{pmatrix}
\Id_k & 0\cr
0 & J_{2n}
\end{pmatrix}$ so that $\sigma=\Ad R$. We have $\sigma \begin{pmatrix}
A & B\cr
B' & C
\end{pmatrix}=\begin{pmatrix}
A & -BJ_{2n}\cr
J_{2n} B' & -J_{2n}CJ_{2n}
\end{pmatrix}$ so that 

\begin{eqnarray*}
\g^\sigma &=& \left\{ \begin{pmatrix}
A & 0\cr
0 & C
\end{pmatrix}\in \g\,\vert\, [C,J_{2n}]=0    \right\}\cr
 \p &=& \left\{ \begin{pmatrix}
0 & B\cr
B' & 0
\end{pmatrix}\in \g \right\} \qquad \qquad \qquad \qquad  J^{\pm}\begin{pmatrix}
0 & B\cr
B' & 0
\end{pmatrix}=\begin{pmatrix}
0 & -BJ_{2n}\cr
J_{2n}B' & 0
\end{pmatrix}\cr
\g^\sigma_{-1} &=& \left\{ \begin{pmatrix}
0 & 0\cr
0 & C
\end{pmatrix}\in \g\,\vert\, CJ_{2n}+J_{2n} C=0    \right\} \qquad J^\pm\begin{pmatrix}
0 & 0\cr
0 & C
\end{pmatrix}= \begin{pmatrix}
0 & 0\cr
0 & \pm J_{2n}C
\end{pmatrix}
\end{eqnarray*}
the last equality following from $\tau =\exp(\frac14 \pi \ad(\rho ))=\Ad \exp(\frac14 \pi \rho)$ so that  
\begin{eqnarray}
\tau \begin{pmatrix}
A & B\cr
B' & C
\end{pmatrix}&=& \Ad\left( \begin{pmatrix}
\Id & 0\cr
0 & \frac{1}{\sqrt{2}}(\Id +J_{2n})
\end{pmatrix} \right)\begin{pmatrix}
A & B\cr
B' & C
\end{pmatrix} \nonumber\\
 &=& \begin{pmatrix}
A & \frac{1}{\sqrt{2}}(B-BJ_{2n})\cr
\frac{1}{\sqrt{2}}(B'+J_{2n}B') & \half(C +J_{2n}C-CJ_{2n}-J_{2n}CJ_{2n})
\end{pmatrix} \label{eq:tau}
\end{eqnarray}
and thus 
$J^\pm\begin{pmatrix}
0 & 0\cr
0 & C
\end{pmatrix}=\pm \tau \begin{pmatrix}
0 & 0\cr
0 & C
\end{pmatrix}= \begin{pmatrix}
0 & 0\cr
0 & \pm J_{2n}C
\end{pmatrix}$ when $CJ_{2n}+J_{2n} C=0$.

\medskip

Let $G$ be the connected Lie subgroup of $Gl(m,\R)$ with Lie algebra $\g$; then $$\tilde\sigma:G\rightarrow G :
g\mapsto RgR^{-1}$$
is an automorphism of $G$ lifting $\sigma$ and $( G,\tilde\sigma,G^{\tilde\sigma}_0)$ is a  $4$-symmetric triple.\\
On the corresponding homogeneous $4$-symmetric space $G/G^{\tilde\sigma}_0$, following Proposition \ref{prop:D}, $J^+$ defines  an  invariant  complex structure $\J^+$ and $J^-$ an invariant almost complex structure $\J^-$ for which the image of the Nijenhuis tensor
is the invariant distribution whose value at the base point is %spanned by
%{\scriptsize{\begin{equation}
%N^{\J^-}_{p_0}\left(\begin{pmatrix}
%0 & B\cr
%B' &  C
%\end{pmatrix},\begin{pmatrix}
%0 & \tilde{B}\cr
%\tilde{B'} &  \tilde{C}
%\end{pmatrix}\right)=\begin{pmatrix}
%0 & -4(B\tilde{C}-\tilde{B}C)\cr
%-4(C\tilde{B'}-\tilde{C}B') &  -2\left((B'\tilde{B}-\tilde{B'}B)+J_{2n}(B'\tilde{B}-\tilde{B'}B)J_{2n})\right)
%\end{pmatrix}.
%\end{equation}}}
%Hence, 
 $\Image N^{\J^-}_{p_0}=[\p,\g^\sigma_{-1}]+[\p,\p]\cap \g^\sigma_{-1}$ and
$\J^-$ is maximally non integrable iff $[\p,\g^\sigma_{-1}]= \p$ and $[\p,\p]\cap \g^\sigma_{-1}=\g^\sigma_{-1}$.
\begin{rem}\label{rmkintegr}
Since {\scriptsize{$\left[\begin{pmatrix}
0 & 0\cr
0 &  C
\end{pmatrix},\begin{pmatrix}
0 & B\cr
B' &  0
\end{pmatrix}\right]=\begin{pmatrix}
0 & -BC\cr
CB' &  0
\end{pmatrix}$}}, the equality $[\p,\g^\sigma_{-1}]= \p$ will automatically hold if there exist $k\ge 1$  elements
{\scriptsize{$\begin{pmatrix}
0 & 0\cr
0 &  \tilde{C}_i
\end{pmatrix}$}} in $\g^\sigma_{-1}$ such that $$\sum_{i=1}^k \tilde{C}_i^2=\Id_{2n}.$$
\end{rem}
The stabilizer of $\rho$ under the   action by conjugation of $G$ in $\gl(m,\R)$ is given by 
\begin{eqnarray*}
\left\{ g\in G\,\vert\, g\rho g^{-1}=\rho \right\}&=&\left\{ g\in G\,\vert\, g\rho=\rho g \right\}\\
&=&\left\{\begin{pmatrix}
E & 0\cr
0 &  D
\end{pmatrix}\in G\,\vert\, DJ_{2n}=J_{2n}D \right\}\\
&=&\left\{ g\in G\,\vert\, gR=Rg \right\}=\left\{ g\in G\,\vert\, \tilde{\sigma}(g)=g \right\},
\end{eqnarray*} 
so coincides with $G^{\tilde\sigma}$. 
Since $J^\pm$ coincides with half the bracket with $\rho$ on $\g^\sigma_{-1}$, it commutes with the adjoint action of any element in $G^{\tilde\sigma}$. 
Hence, on the  homogeneous $4$-symmetric space $G/G^{\tilde\sigma}$, corresponding to the  $4$-symmetric triple $( G,\tilde\sigma,G^{\tilde\sigma})$,  $J^+$  defines also an  invariant  complex structures $\J^+$ and $J^-$ also an invariant almost complex structure $\J^-$ for which the image of the Nijenhuis tensor
is the invariant distribution whose value at the base point is given as above.

%%%%%%%%%%%%%%%%%%%%%%%%%%%%%%%%%

\subsubsection{The symplectic structure}
Since the stabilizer of $\rho$ under the   action by conjugation of $G$ in $\gl(m,\R)$ coincides with $G^{\tilde\sigma}$, the homogeneous space  $G/G^{\tilde\sigma}$
is  diffeomorphic to the  orbit of   $\rho$ under the action of $G$  in $\gl(m,\R)$. 

Let us define the $G$-invariant $2$-form $\omega$ whose value at the base point is given by 
\begin{equation}\label{eq:omega}
\tilde{\Omega}(X,Y):= \Tr(\rho[X,Y]) \qquad \forall X,Y \in \g^\sigma_{-1}+\p.
\end{equation}
Observe that $\Tr(\rho[X,Y])=0$ for any $X\in \g^\sigma$ and any $Y\in \g$.
The $2$-form is  invariant by the symmetries since 
$$
\tilde{\Omega}(\sigma X,\sigma Y)= \Tr(\rho \Ad R ([X,Y]))=\Tr( \Ad R (\rho [X,Y]))=\Tr(\rho  [X,Y])=\tilde{\Omega}( X, Y).
$$
The $2$-form $\omega$ is closed since, for all $X,Y,Z \in \g^\sigma_{-1}+\p$, denoting by $pr:\g\rightarrow \g^\sigma_{-1}+\p$ the projection parallel to $\g^\sigma$, we have: 
$$
\cyclic_{XYZ} \tilde{\Omega}([X,Y]_{pr},Z)=\cyclic_{XYZ} \Tr(\rho[[X,Y]_{pr},Z])=\cyclic_{XYZ} \Tr(\rho[[X,Y],Z])=0.
$$
The $2$-form $\omega$ is an invariant symplectic structure on the $4$-symmetric space defined by $( G,\tilde\sigma,G^{\tilde\sigma})$ and on its cover defined by $( G,\tilde\sigma,G^{\tilde\sigma}_0)$, if and only if it is non degenerate.
This will be true if and only if
$$
\tilde{\Omega}':\p\times \p\rightarrow \R  \quad \textrm{ with }\quad 
\tilde{\Omega}'(\begin{pmatrix}
0 & B_1\cr
B'_1 &  0
\end{pmatrix},\begin{pmatrix}
0 & B_2\cr
B'_2 &  0
\end{pmatrix})=\Tr \left(J_{2n}(B'_1B_2-B'_2B_1)\right)
$$
and 
$$
\tilde{\Omega}^v:\g^\sigma_{-1}\times \g^\sigma_{-1}\rightarrow \R \quad \textrm{ with }\quad 
\tilde{\Omega}^v(\begin{pmatrix}
0 & 0\cr
0 &  C_1
\end{pmatrix},\begin{pmatrix}
0 & 0\cr
0 &  C_2
\end{pmatrix})=2\Tr \left(J_{2n}C_1C_2\right)
$$ 
are non degenerate, i.e. iff
$$
g' :\p\times \p\rightarrow \R :X,Y \mapsto \Tr(XY) \quad \textrm{and}\quad  \beta^v : \g^\sigma_{-1} \times \g^\sigma_{-1} \rightarrow \R : X,Y \mapsto \Tr(XY)
$$ are non degenerate. This we shall now assume.\\
If $\rho$ belongs to the Lie algebra  $\g$ and if the map
$ \beta :\g \times \g \rightarrow \R : X,Y \mapsto \Tr(XY)$ is non degenerate (which is the case if $\g$ is simple), then  $G/G^{\tilde{\sigma}}$,
identifies with the coadjoint orbit of the element $\rho^\flat: \g \rightarrow \R : X \mapsto \rho^\flat:=\Tr \rho X$ and the $2$-form $\omega$ is automatically non degenerate because it is the Kirillov-Kostant-Souriau symplectic $2$-form on this orbit.\\

In all cases, the natural almost complex structures $\J^\pm$ are compatible with $\omega$ since
\begin{eqnarray*}
\tilde{\Omega}(\begin{pmatrix}
0 & B_1\cr
B'_1 &  C_1
\end{pmatrix},\begin{pmatrix}
0 & B_2\cr
B'_2 &  C_2
\end{pmatrix})&=&\Tr \left(J_{2n}( 2C_1C_2+B'_1B_2-B'_2B_1)\right)\\
 &=&\tilde{\Omega}(\begin{pmatrix}
0 & -B_1 J_{2n}\cr
J_{2n} B'_1 & \pm J_{2n} C_1
\end{pmatrix},\begin{pmatrix}
0 & -B_2J_{2n}\cr
J_{2n} B'_2 & \pm J_{2n} C_2
\end{pmatrix}).
\end{eqnarray*}
In general, they are neither positive nor negative :
\begin{eqnarray*}
\tilde{\Omega}(\begin{pmatrix}
0 & B\cr
B'&  C
\end{pmatrix},J^\pm\begin{pmatrix}
0 & B\cr
B' &  C
\end{pmatrix})&=& \Tr \begin{pmatrix}
0 & 0\cr
0&  J_{2n}
\end{pmatrix}\left[\begin{pmatrix}
0 & B\cr
B'&  C
\end{pmatrix},\begin{pmatrix}
0 & -BJ_{2n}\cr
J_{2n}B'&  \pm J_{2n}C
\end{pmatrix}\right]\\
&=&
   \Tr\left(J_{2n}(-B'BJ_{2n}\pm CJ_{2n}C-J_{2n}B'B\mp J_{2n}C^2)\right)\\
    &=&2\Tr B'B\pm 2\Tr C^2,
\end{eqnarray*}
 and $\pm 2\Tr C^2=\pm 4(\Tr (c^2)+\Tr(d^2))$ for $C=\begin{pmatrix}
c & d\cr
d &  -c
\end{pmatrix}$.

%%%%%%%%%%%%%%%%%%%%%%%%%
\subsubsection{The Levi Civita  connection and  the Ricci  Hermitian property }

We consider as above the homogeneous space $G/G^{\tilde{\sigma}}$ or $G/G^{\tilde{\sigma}}_0$
and we assume that the  $2$-form  $\tilde{\Omega}$ given in equation \ref{eq:omega} is non degenerate.
The pseudo-Riemannian metrics associated to $\J^\pm $ will be denoted by $g^\pm $; they are the 
$G$-invariant metrics whose value at the base point are given by
$$
G^\pm(X,Y):=\tilde{\Omega}(X,J^\pm Y):= \Tr(\rho[X,J^\pm Y ]) \qquad \forall X,Y \in \g^\sigma_{-1}+\p,
$$
thus
\begin{eqnarray*}
G^\pm(\begin{pmatrix}
0 & B_1\cr
B'_1&  C_1
\end{pmatrix},\begin{pmatrix}
0 & B_2\cr
B'_2 &  C_2
\end{pmatrix})&=& \Tr \begin{pmatrix}
0 & 0\cr
0&  J_{2n}
\end{pmatrix}\left[\begin{pmatrix}
0 & B_1\cr
B'_1&  C_1
\end{pmatrix},\begin{pmatrix}
0 & -B_2J_{2n}\cr
J_{2n}B'_2&  \pm J_{2n}C_2
\end{pmatrix}\right]\\
%&=& \Tr\left(J_{2n}(-B'_1B_2J_{2n}\pm C_1J_{2n}C_2-J_{2n}B'_2B_1\mp J_{2n}C_2C_1)\right)\\
    &=&\Tr (B'_1B_2+B'_2B_1)\pm 2\Tr C_1C_2.
\end{eqnarray*}
The corresponding Levi Civita connection are denoted by $\nabla^{g^\pm}$. Invariance implies, $(\mathcal{L}_{A^*}g^\pm)(B^*,C^*)={A^*}g^\pm(B^*,C^*)-g^\pm([A^*,B^*],C^*)-g^\pm(B^*,[A^*,C^*])=0$, and  
$$
2g^\pm(\nabla^{g^\pm}_{A^*}B^*,C^*)=g^\pm([A^*,B^*],C^*)+g^\pm([A^*,C^*],B^*)+g^\pm([B^*,C^*],A^*),
$$
for any $A,B,C \in \g$, where $A^*$ denotes the fundamental vector field associated to $A\in \g$, i.e.    
 $A^*_p:=\frac{d}{dt} \exp -tA.p\vert_{t=0}$ and,  at the base point $p_0=\pi(e)$,  $A^*_{p_0}=-\pi_* (A)$.
 In our case, it gives, at the base point, for all $X_i$'s in $\g^\sigma_{-1}$ and all $Y_i$'s in $\p=\g^{\sigma^2}_{-1}$:
 $$
 \nabla^{g^\pm}_{X_1^*}X_2^*(p_0)= 0\qquad \nabla^{g^\pm}_{X^*}Y^*(p_0)=\pm [X,Y]^*_{p_0} \qquad \nabla^{g^\pm}_{Y_1^*}Y_2^*(p_0)= \half [Y_1,Y_2]^*_{p_0}
 $$
 and $\nabla^{g^\pm}$ being torsion free, $\nabla^{g^+}_{Y^*}X^*(p_0)=0$ whereas $ \nabla^{g^-}_{Y^*}X^*(p_0)=-2[X,Y]^*_{p_0}$.
 
%%%%%%%%%%%%%%%%%%%%%%%%%%%%%%%%
Since $\J^+$ is integrable, we know that $\nabla^{g^+}\J^+=0$ hence $\J^+$ commutes with the curvature
$R^{g^+}(X,Y)=\nabla^{g^+}_X\circ \nabla^{g^+}_Y-\nabla^{g^+}_Y\circ \nabla^{g^+}_X-\nabla^{g^+}_{[X,Y]}$.\\ Since $g^+(R^{g^+}(X,Y)Z,T)=g^+(R^{g^+}(Z,T)X,Y)$, this implies that
 $$R^{g^+}(\J^+X,\J^+Y)=R^{g^+}(X,Y).$$
  The Ricci tensor $Ric^{g^+}(X,Z)=\Tr [Y\rightarrow R^{g^+}(X,Y)Z]$ is thus  hermitian for $\J^+$ :
$$
Ric^{g^+}(\J^+X,\J^+Z)=Ric^{g^+}(X,Z) \qquad \forall X,Z.
$$ 
A direct computation at the base point shows that we have a similar result for $\J^-$. One uses the fact that the curvature of a torsion free connection in a homogeneous situation can be computed from
$$
R^{g^\pm}({A^*},B^*)C^*=\nabla^{g^\pm}_{\nabla^{g^\pm}_{C^*}B^*}A^*-\nabla^{g^\pm}_{\nabla^{g^\pm}_{C^*}A^*}B^*+\nabla^{g^\pm}_{[A,B]^*}C^*-[[A,B],C]^*.
$$
 At the base point, for all $X_i$'s in $\g^\sigma_{-1}$ and all $Y_i$'s in $\p=\g^{\sigma^2}_{-1}$, it gives:
 {\scriptsize{
\begin{eqnarray*}
R_{p_0}^{g^\pm}({X_1^*},X_2^*)X_3^*&=&-[[X_1,X_2],X_3]^*_{p_0}   \\
R_{p_0}^{g^+}({X_1^*},X_2^*)Y^*&=&-[[X_1,X_2],Y]^*_{p_0} \\ 
R_{p_0}^{g^-}({X_1^*},X_2^*)Y^*&=&3[[X_1,X_2],Y]^*_{p_0} \\ 
R_{p_0}^{g^+}({X_1^*},Y^*)X_2^*&=&-[[X_1,Y],X_2]^*_{p_0} \\ 
R_{p_0}^{g^-}({X_1^*},Y^*)X_2^*&=&(-[[X_1,Y],X_2]+2[[X_1,X_2]Y])^*_{p_0} \\ 
R_{p_0}^{g^+}({X^*},Y_1^*)Y_2^*&=&-\half[[X,Y_1],Y_2]^*_{p_0} \\ 
R_{p_0}^{g^-}({X^*},Y_1^*)Y_2^*&=&(-\half[[X,Y_1],Y_2]+[[X,Y_2]Y_1])^*_{p_0} \\ 
R_{p_0}^{g^+}(Y_1^*,Y_2^*)X^*&=&-\half[[Y_1,Y_2],X]^*_{p_0} \\ 
R_{p_0}^{g^-}(Y_1^*,Y_2^*)X^*&=&-{{\frac{3}{2}}}[[Y_1,Y_2],X]^*_{p_0} \\ 
R_{p_0}^{g^+}({Y_1^*},Y_2^*)Y_3^*&=&\frac{1}{4}(-[[Y_1,Y_2],Y_3]+[\sigma[Y_3,Y_1],Y_2]+[\sigma[Y_2,Y_3],Y_1]-2[\sigma[Y_1,Y_2],Y_3] )^*_{p_0}   \\
R_{p_0}^{g^-}({Y_1^*},Y_2^*)Y_3^*&=&\frac{1}{4}(-7[[Y_1,Y_2],Y_3]-[\sigma[Y_3,Y_1],Y_2]-[\sigma[Y_2,Y_3],Y_1]+2[\sigma[Y_1,Y_2],Y_3] )^*_{p_0}.
\end{eqnarray*}}}
The Ricci tensor is invariant and its value at the base point is given by
 {\scriptsize{
 \begin{eqnarray*}
Ric_{p_0}^{g^+}(X_1^*,X_2^*)&=&\Tr\vert_{\g^\sigma_{-1}} \ad X_2\circ \ad X_1+\Tr\vert_{\p} \ad X_2\circ \ad X_1  \\
Ric_{p_0}^{g^+}(X^*,Y^*)&=&0  \\
Ric_{p_0}^{g^+}(Y_1^*,Y_2^*)&=&\half\Tr\vert_{\g^\sigma_{-1}}\half(\id-\sigma)\circ  \ad Y_2\circ \ad Y_1+\\
&&\frac{1}{4}\Tr\vert_{\p}( \ad Y_2\circ \ad Y_1 +\ad(\sigma[Y_2,Y_1])+\ad Y_1\circ \sigma\circ \ad Y_2 +2\ad Y_2\circ \sigma\circ \ad Y_1  )\\
Ric_{p_0}^{g^-}(X_1^*,X_2^*)&=&\Tr\vert_{\g^\sigma_{-1}} \ad X_2\circ \ad X_1+\Tr\vert_{\p} \ad X_2\circ \ad X_1+2 \Tr\vert_{\p} \ad[X_1,X_2] \\
&=& \Tr\vert_{\g^\sigma_{-1}} \ad X_2\circ \ad X_1+\Tr\vert_{\p} \ad X_2\circ \ad X_1\\
Ric_{p_0}^{g^-}(X^*,Y^*)&=&0 \\
Ric_{p_0}^{g^-}(Y_1^*,Y_2^*)&=&\Tr\vert_{\g^\sigma_{-1}}\half(\id-\sigma)( \half\ad Y_2\circ \ad Y_1-\ad Y_1\circ \ad Y_2)\\
&&+\frac{1}{4}\Tr\vert_{\p}( 7\ad Y_2\circ \ad Y_1-\ad(\sigma[Y_2,Y_1])-\ad Y_1\circ \sigma\circ \ad Y_2-2\ad Y_2\circ \sigma\circ \ad Y_1 ). \\
\end{eqnarray*} }}
The Ricci tensor $Ric_{p_0}^{g^+}$ is  $\J^+$ Hermitian and the tensor $Ric_{p_0}^{g^-}$ is  $\J^-$ Hermitian because each term is invariant under $J^\pm$ which coincide with $\sigma$ on $\p$ and with $\pm \tau=\pm \exp\frac{\pi}{4}D$ on $\g^\sigma_{-1} $.
Indeed, since $\ad \sigma Y= \sigma\circ \ad Y\circ \sigma^{-1}$, on any subspace $V$ stable by $\sigma$ (in particular on $\p$ and on $\g^\sigma_{-1}$), one has
 {\scriptsize{
 \begin{eqnarray*}
\Tr\vert_{V} \ad \sigma Y_2\circ \ad \sigma Y_1 &=& \Tr\vert_{V}\sigma\circ\ad  Y_2\circ \ad  Y_1\circ \sigma^{-1}=\Tr\vert_{V}\ad  Y_2\circ \ad  Y_1\\
\Tr\vert_{V} \ad \sigma Y_2\circ\sigma \circ \ad \sigma Y_1 &=& \Tr\vert_{V}\sigma\circ\ad  Y_2\circ \sigma\circ \ad  Y_1\circ \sigma^{-1}=\Tr\vert_{V}\ad  Y_2\circ\sigma \circ \ad  Y_1\\
\end{eqnarray*} }}
so that $Ric_{p_0}^{g^\pm}(\J^\pm Y_1^*,\J^\pm Y_2^*)=Ric_{p_0}^{g^\pm}(Y_1^*,Y_2^*)$,
and, similarly, since $\ad \tau X= \tau\circ \ad X\circ \tau^{-1}$, on any subspace $V$ stable by $\tau$ (in particular on $\p$ and on $\g^\sigma_{-1}$), one has
{\scriptsize{
 \begin{eqnarray*}
\Tr\vert_{V} \ad \tau X_2\circ \ad \tau X_1 &=& \Tr\vert_{V}\tau\circ\ad  X_2\circ \ad  X_1\circ \tau^{-1}=\Tr\vert_{V}\ad  X_2\circ \ad  X_1\\
\end{eqnarray*} }}
so that $Ric_{p_0}^{g^\pm}(\J^\pm X_1^*,\J^\pm X_2^*)=Ric_{p_0}^{g^\pm}(X_1^*,X_2^*)$.

\medskip
For elements $X_i$'s in $\g_{-1}^\sigma$ of the form $X_i=\left(\begin{smallmatrix} 0 & 0\cr  0&  C_i  \end{smallmatrix}\right)$ and $Y_i$'s in $\p$ of the form $Y_i=\left(\begin{smallmatrix}  0 & B_i\cr B_i'&  0\end{smallmatrix}\right)$, one has
{\scriptsize{ \begin{eqnarray}\label{eq:Einstein}
Ric_{p_0}^{g^+}(X_1^*,X_2^*)&=&\Tr\vert_{\g^\sigma_{-1}} \ad X_2\circ \ad X_1+\Tr\vert_{\p} \ad X_2\circ \ad X_1  \nonumber\\
&=&\Tr\left( C\rightarrow (C_2C_1C+CC_1C_2-C_1CC_2-C_2CC_1)\right)+\Tr \left( (B,B')\rightarrow (BC_1C_2, C_2C_1B')\right)\label{eq:Einstein1}\\
Ric_{p_0}^{g^+}(X_1^*,Y_1^*)&=&0  \quad \textrm{and} \quad  g^+_{p_0}(X_1^*,Y_1^*)=0;\label{eq:Einstein2}\\
Ric_{p_0}^{g^+}(Y_1^*,Y_2^*)&=&\half\Tr\vert_{\g^\sigma_{-1}} \half(\id-\sigma)\ad Y_2\circ \ad Y_1
+\frac{1}{4}\Tr\vert_{\p}( \ad Y_2\circ \ad Y_1 +\ad(\sigma[Y_2,Y_1]))\nonumber\\
&&+\frac{1}{4}\Tr\vert_{\p}( \ad Y_1\circ \sigma\circ \ad Y_2 +2\ad Y_2\circ \sigma\circ \ad Y_1  )\nonumber\\
&=&\half\Tr\left( C\rightarrow \half((B'_2B_1-J_{2n}B'_2B_1J_{2n})C+C(B'_1B_2-J_{2n}B'_1B_2J_{2n})\right)\label{eq:Einstein3}\\
&&+\frac{1}{4}\Tr\Bigg( (B,B')\rightarrow \Big( (2B_2B'_1-B_1B'_2)B+B(3B'_1B_2+B'_2B_1-JB'_1B_2J +JB'_2B_1J)\nonumber\\
&&\qquad\qquad \qquad  -(B_1JB'_2+2B_2JB'_1)BJ-2B_2B'B_1-3B_1B'B_2+B_1JB'B_2J+2 B_2JB'B_1J  \, ,\nonumber\\
&&\qquad\qquad \qquad  B'(2B_1B'_2-B_2B'_1)+(3B'_2B_1+B'_1B_2+JB'_1B_2J -JB'_2B_1J) B'      \nonumber   \\
&&\qquad\qquad \qquad  -JB'(B_2JB'_1+2B_1JB'_2)-2B'_1BB'_2-3B'_2BB'_1+JB'_2BJB'_1+2 JB'_1BJB'_2 \Big)\Bigg)\nonumber
 \end{eqnarray}}}
where the traces are computed on elements  $C\in Mat(2n\times 2n)$ of the form $C=\left(\begin{smallmatrix}  c & c'\cr c'&  -c\end{smallmatrix}\right)$
such that $\left(\begin{smallmatrix} 0 & 0\cr  0&  C  \end{smallmatrix}\right) \in \g_{-1}^\sigma $, 
and on elements $B \in Mat(k\times 2n),B'\in Mat(2n\times k)$ 
such that  $\left(\begin{smallmatrix}  0 & B\cr B'&  0\end{smallmatrix}\right) \in \p$. Observe that,
for any $D_i=\left(\begin{smallmatrix}  d_i & d'_i\cr -d'_i&  d_i\end{smallmatrix}\right)$,
{\scriptsize{
\begin{eqnarray}
\Tr\Big( C\rightarrow (D_1C+CD_2)\Big)  &=& \Tr\Big( (c,c')\rightarrow (d_1c+d'_1c'+cd_2-c'd'_2,d_1c'-d'_1c+cd'_2+c'd_2)\Big) \\
\Tr\Big( C\rightarrow C_1CC_2\Big)&=&\Tr\Big( (c,c')\rightarrow (c_1cc_2+c'_1c'c_2+c_1c'c'_2-c'_1cc'_2, c'_1cc_2-c_1c'c_2+c_1cc'_2+c'_1c'c'_2 )\Big)\nonumber\\
\end{eqnarray}
}}
We shall use those identities to show that $g^+$ is Einstein in the examples given in section \ref{section:examples}.

%%%%%%%%%%%%%%%%%%%%%%%%%
\subsubsection{The Chern  connection, the  Chern  Ricci form and the property of being special }
The almost complex structure $\J^+$ is integrable, so $\nabla^{g^+}\J^+=0$, and the Chern connection $\nabla^{C^+}$ in that case coincides with the Levi Civita connection $\nabla^{g^+}$.
The Chern Ricci form  \eqref{eq:ChernRicci} is thus given by
$$
ChernRicci^{\J^+}(X,Y):= \Tr \J^+ R^{g^+}(X,Y) = Ric^{g^+}(X,\J^+Y)-Ric^{g^+}(Y,\J^+X).
$$
Since we have the Ricci Hermitian property, this yields
$$
ChernRicci^{\J^+}(X,Y):= 2 Ric^{g^+}(X,\J^+Y).
$$
And $\J^+$ is special (i.e. $ChernRicci^{\J^+}$ is proportional to $\omega$) if and only if
$g^+$  is Einstein, i.e. $Ric^{g^+}$ is proportional to ${g^+}$.  

\medskip

For the almost complex structure $\J^-$,  the Chern connection $\nabla^{C^-}$ is given by
$$
\nabla^{C^-}_XY=\nabla^{g^-}_XY-\half \J (\nabla^{g^-}\J)_XY=\half \nabla^{g^-}_XY-\half \J \nabla^{g^-}_X(\J Y) \qquad \forall X,Y \in \mathfrak{X}(M);
$$
it is invariant and given, for all $A,B \in \g$  by
$$
\nabla^{C^-}_{A^*}B^*=\half \nabla^{g^-}_{A^*}B^*-\half \J^-(\nabla^{g^-}_{\J^- B^*}A^*+[A^*,\J^- B^*])=\half \nabla^{g^-}_{A^*}B^* +\half [{A^*},B^*]-\half \J^-\nabla^{g^-}_{\J^- B^*}A^*.
$$
At the base point $p_0$, for all $X_i$'s in $\g^\sigma_{-1}$ and all $Y_i$'s in $\p=\g^{\sigma^2}_{-1}$,  the Chern connection is given, by
\begin{eqnarray*}
\nabla^{C^-}_{X_i^*}X_k^*(p_0)&=& 0, \\
\nabla^{C^-}_{X_i^*}Y_k^*(p_0)&=& \left(-\half[X_i,Y_k]+\half[X_i,Y_k]+ \sigma [X_i,\sigma Y_k]\right)^*_{p_0}=[ \sigma X_i,-Y_k]^*_{p_0}= [X_i,Y_k]^*_{p_0}\\
\nabla^{C^-}_{Y_k^*}X_i^*(p_0)&=& \left(\frac{3}{2}[Y_k,X_i]-\half \sigma [\tau X_i,Y_k,]\right)^*_{p_0}= [Y_k,X_i]^*_{p_0}, \\
\nabla^{C^-}_{Y_i^*}Y_k^*(p_0)&=& \left(\frac{3}{4}[Y_i,Y_k]^*_{p_0}+\frac{1}{4} \tau
\frac{1-\sigma}{2}[{\sigma Y_k},Y_i] \right)=[Y_i,Y_k]^*_{p_0}.
\end{eqnarray*}
because,  $\sigma[\tau X_i,Y_k]=- [X_i,Y_k] $ and $\tau (\Id-\sigma)[{\sigma Y_k},Y_i]=(\Id-\sigma)[ Y_k,Y_i]$.\\
Hence we have, for all $A,B \in \g^\sigma_{-1}+\p$
\begin{equation}
\nabla^{C^-}_{A^*}B^*(p_0)=[A,B]^*_{p_0}\quad \textrm{ and  }\quad T^{\nabla^{C^-}}_{p_0}({A^*},B^*)=[A,B]^*_{p_0}=\frac{1}{4} N^{\J^-}_{p_0}({A^*},B^*).
\end{equation}
The curvature of a connection with torsion in a homogeneous situation can be computed from
\begin{eqnarray*}
R^{C^-}({A^*},B^*)C^*&=&\nabla^{{C^-}}_{\nabla^{{C^-}}_{B^*}C^*}A^*-\nabla^{{C^-}}_{\nabla^{{C^-}}_{A^*}C^*}B^* +{\nabla^{{C^-}}_{[A,B]^*}C^*}+{\nabla^{{C^-}}_{[A,C]^*}{B^*}}-{\nabla^{{C^-}}_{[B,C]^*}{A^*}}\\
&&-[[A,B],C]^*+T^{\nabla^{C^-}}({A^*},{\nabla^{{C^-}}_{B^*}C^*})-T^{\nabla^{C^-}}({B^*},{\nabla^{{C^-}}_{A^*}C^*}) \\
&& +T^{\nabla^{C^-}}(B^*,[A,C]^*)-T^{\nabla^{C^-}}(A^*,[B,C]^*).
\end{eqnarray*}
At the base point $p_0$, for all $A,B,C \in \g^\sigma_{-1}+\p$, we have 
\begin{eqnarray*}
R^{C^-}_{p_0}({A^*},B^*)C^*&=&\nabla^{{C^-}}_{[{B},C]^*}A^*({p_0})-\nabla^{{C^-}}_{[{A},C]^*}B^* ({p_0})+{\nabla^{{C^-}}_{[A,B]^*}C^*}({p_0})+{\nabla^{{C^-}}_{[A,C]^*}{B^*}}({p_0})\\ &&-{\nabla^{{C^-}}_{[B,C]^*}{A^*}}({p_0})
-[[A,B],C]^*_{p_0}+T^{\nabla^{C^-}}_{p_0}({A^*},[B,C]^*)\\ &&-T^{\nabla^{C^-}}_{p_0}(B^*,[A,C]^*) 
 +T^{\nabla^{C^-}}_{p_0}(B,^*[A,C]^*)-T^{\nabla^{C^-}}_{p_0}(A^*,[B,C]]^*)\\
 &=& {\nabla^{{C^-}}_{[A,B]^*}C^*}({p_0})-[[A,B],C]^*_{p_0}\\
&=&  \left( [ pr[A,B],C]-[[A,B],C] \right)^*_{p_0}
\end{eqnarray*}
where $pr:\g\rightarrow \g^\sigma_{-1}+\p$ is the projection parallel to $\g^\sigma$.
Hence, for all $X_i$'s in $\g^\sigma_{-1}$ and all $Y_i$'s in $\p=\g^{\sigma^2}_{-1}$, since
$pr[X_i,X_j]=0, pr[X_i,Y_j]=[X_i,Y_j]$ and $pr[Y_i,Y_j]=\half(\Id -\sigma)[Y_i,Y_j]$, the Chern curvature for $\J^-$ is given by
{\scriptsize{\begin{eqnarray*}
R^{C^-}_{p_0}({X_i^*},X_j^*)X_k^*&=&-([[X_i,X_j],X_k])^*_{p_0}\\
R^{C^-}_{p_0}({X_i^*},X_j^*)Y_k^*&=&-([[X_i,X_j],Y_k])^*_{p_0}\\
R^{C^-}_{p_0}({X_i^*},Y_k^*)X_j^*&=&0\\
R^{C^-}_{p_0}({X_i^*},Y_j^*)Y_k^*&=&0\\
R^{C^-}_{p_0}({Y_i^*},Y_j^*)X_k^*&=&-(\half[(\Id +\sigma)[Y_i,Y_j],X_k])^*_{p_0}\\
R^{C^-}_{p_0}({Y_i^*},Y_j^*)Y_k^*&=&-(\half[(\Id +\sigma)[Y_i,Y_j],Y_k])^*_{p_0}.
\end{eqnarray*}}}
Hence the  Chern Ricci form, $ChernRicci^{\J^-}(X,Y):= \Tr \J^- R^{C^-}(X,Y)$,  for elements $X_i$'s in $\g_{-1}^\sigma$ of the form $X_i=\left(\begin{smallmatrix}
0 & 0\cr
0&  C_i
\end{smallmatrix}\right)$ and $Y_i$'s in $\p$ of the form $Y_i=\left(\begin{smallmatrix}
0 & B_i\cr
B_i'&  0
\end{smallmatrix}\right)$, reads
{\scriptsize{\begin{eqnarray}
ChernRicci^{\J^-}_{p_0}({X_i^*},X_j^*)&=&\Tr\vert_{\g_{-1}^\sigma} \tau \circ \ad [X_i,X_j]-\Tr\vert_{\p} \,\sigma \circ \ad [X_i,X_j]\nonumber\\
     &=&\Tr\left( C\rightarrow D(X_i,X_j)C+C D(X_i,X_j)\right)\nonumber\\
     && -\Tr\left((B,B')\rightarrow (BD(X_i,X_j), D(X_i,X_j)B')\right)\label{eq:RicChern1}\\
     \textrm{with}\, D(X_i,X_j)&=&J_{2n}[C_i,C_j]\nonumber\\
ChernRicci^{\J^-}_{p_0}({X_i^*},Y_k^*)&=&0 \label{eq:RicChern2}\\
ChernRicci^{\J^-}_{p_0}({Y_i^*},Y_j^*)&=&\Tr\vert_{\g_{-1}^\sigma} \tau \circ \ad (\frac{\Id+\sigma}{2}[Y_i,Y_j])-\Tr\vert_{\p}\, \sigma \circ \ad (\frac{\Id+\sigma}{2}[Y_i,Y_j])\nonumber\\
&=&\Tr\left( C\rightarrow D(Y_i,Y_j)C+CD(Y_i,Y_j)\right)\nonumber\\
&& -\Tr\left((B,B')\rightarrow (BD(Y_i,Y_j), D(Y_i,Y_j)B')\right)\label{eq:RicChern3}\\
 \textrm{with}\,D(Y_i,Y_j) &=&\half( J_{2n}(B'_iB_j-B'_jB_i)+(B'_iB_j-B'_jB_i)J_{2n})\nonumber
\end{eqnarray}}}
for $C$ such that $\left(\begin{smallmatrix}
0 & 0\cr
0&  C
\end{smallmatrix}\right) \in \g_{-1}^\sigma$  and $(B,B')$ such that $\left(\begin{smallmatrix}
0 & B\cr
B'&  0
\end{smallmatrix}\right)\in \p$.
Remark that $C$ is of the form $\left(\begin{smallmatrix}
c & c'\cr
c' & -c
\end{smallmatrix}\right)$ and $D$   is of the form $\left(\begin{smallmatrix}
d & d'\cr
-d' & d
\end{smallmatrix}\right)$.
\medskip

If one can find a basis $\{e_k\}$ of the $2n\times 2n$ matrices $C$ such that $\left(\begin{smallmatrix}
0 & 0\cr
0&  C
\end{smallmatrix}\right) \in \g_{-1}^\sigma$ and of the pairs of $k\times 2n$ and $2n\times k$ matrices $(B,B')$ such that $\left(\begin{smallmatrix}
0 & B\cr
B'&  0
\end{smallmatrix}\right)\in \p$, for which each element is a  linear combination of elementary matrices 
$E_{ij}$ (the  matrix whose only non vanishing entry is a $1$ at the intersection of the $i^{th}$ row and the $j^{th}$ column) with different indices $i$ of lines and different indices $j$ of columns, then, since
the only coefficient of $E_{ij}$ in $De_k$ with $e_k=E_{ij}+\ldots$ is $D_{ii}$ and in $e_kD$ is $D_{jj}$, we obtain that  $\Tr\left( C\rightarrow DC+CD\right)$
 and $\Tr\left((B,B')\rightarrow (BD, DB')\right)$ are two multiples of $\Tr D$.
 Since 
 \begin{eqnarray*}
 \Tr D(X_i,X_j)&=&\Tr J_{2n}[C_i,C_j]=\Tr \rho[X_i,X_j]=\tilde\Omega(X_i,X_j)\\
 \Tr D(Y_i,Y_j) &=&\Tr  J_{2n}(B'_iB_j-B'_jB_i)=\Tr \rho[Y_i,Y_j]=\tilde\Omega(Y_i,Y_j),
 \end{eqnarray*}
 we conclude that  $ChernRicci^{\J^-}$ is proportional to $\omega$, hence $\J^-$ is special.
 This will be done explicitely in the  examples given in the next section.
%%%%%%%%%%%%%%%%%%%%%%%%%%%%%%%%%%
%%%%%%%%%%%%%%%%%%%%%%%%%%%%%%%%%
\section{Examples}\label{section:examples}

\subsection{ $G=Gl(k+2n,\R)$ or $Sl(k+2n,\R)$}

For the group $G=Gl(k+2n,\R)$ (resp. $G=Sl(k+2n,\R)$), $\rho$ lies in its Lie algebra and the map
$\beta: \g \times \g \rightarrow \R : X,Y \mapsto \Tr(XY)$ is non degenerate. In both cases, the $4$-symmetric space
$M$ associated to the triple $(G,\tilde{\sigma},G^{\tilde{\sigma}})$,
$$
G/G^{\tilde{\sigma}}=Gl(k+2n,\R)/\left(Gl(k,\R)\times Gl(n,\mathbb{C})\right)
=Sl(k+2n,\R)/S(Gl(k,\R)\times Gl(n,\mathbb{C})),
$$ 
identifies with the coadjoint orbit of the element $\rho^\flat: \g \rightarrow \R : X \mapsto \Tr \rho X$. The $2$-form $\omega$ is the Kirillov-Kostant-Souriau symplectic $2$-form on this orbit.\\
With the pair $(\omega,\J^+)$, this orbit  is an invariant pseudo-K\"ahler  $4$-symmetric space $M$.\\
With the pair $(\omega,\J^-)$ it is an invariant almost-pseudo-K\"ahler $4$-symmetric space,
and  $\J^-$ is maximally non integrable, since  $[\p,\p]\cap \g^\sigma_{-1} =\g^\sigma_{-1}$ and $[\p,\g^\sigma_{-1}]= \p$. \\
These  equalities follow from remarks \ref{rmkintegr1} and  \ref{rmkintegr}, the fact that 
$\mathfrak{sl}(k+2n,\R)$ is simple, and the fact that 
$\tilde C=\left(\begin{smallmatrix}
0&\Id_n\\\Id_n& 0
\end{smallmatrix}\right)$ satisfies $\tilde C^2=\Id_{2n}$ and $\left(\begin{smallmatrix}
0&0\\0&\tilde C
\end{smallmatrix}\right)$ is in $\g^{\sigma}_{-1}.$

Since $\g^{\sigma}_{-1}=\{ \left(\begin{smallmatrix}
0 & 0\cr
0&  C
\end{smallmatrix}\right)\,\vert \,CJ_{2n}+J_{2n} C= 0 \}$, a basis of the corresponding $C's$ is given by 
$\{E_{i\, j}-E_{n+i\, n+j}, E_{n+i \, j}+E_{i\, n+j}, 1\le i,j\le n\}$ so that, for any $D=\left(\begin{smallmatrix}
d & d'\cr
-d' & d
\end{smallmatrix}\right)$,  
{\scriptsize{$$\Tr\left( C\rightarrow CD\right)=\sum_{1\le i,j\le n}(d_{jj}+d_{jj})=2n \Tr d= n \Tr D,\qquad
\Tr\left( C\rightarrow DC\right)=\sum_{1\le i,j\le n}(d_{ii}+d_{ii})=n \Tr D.$$}}
Since $\p=\{ \left(\begin{smallmatrix}
0 & B\cr
B'&  0
\end{smallmatrix}\right)\}$, a basis of the corresponding set of pairs $(B,B')$ is given by 
$\{(E_{r\, \ell}, E_{\ell'\,  s}), 1\le \ell, \ell' \le 2n , 1\le r,s \le k \}$ and , for any $F_1,F_2\in Mat(2n\times 2n), G_1,G_2 \in Mat(k\times k)$ 
{\scriptsize{
\begin{eqnarray*}
\Tr\left((B,B')\rightarrow (G_1BF_1, F_2B'G_2)\right)
&=&\sum_ {1\le \ell \le 2n , 1\le r \le k}
\left((G_1)_{rr}(F_1)_{\ell \ell}+(F_2)_{\ell \ell}(G_2)_{rr}\right)\\
&=&\Tr G_1 \Tr F_1+\Tr F_2\Tr G_2.
\end{eqnarray*}
}}
The metric $g^+$ is Einstein : $Ric^{g^+}=(k+n)\, g^+$. Indeed, using the identities \eqref{eq:Einstein1}, \eqref{eq:Einstein2} and  \eqref{eq:Einstein3}, one has :
{\scriptsize{\begin{eqnarray*}
Ric_{p_0}^{g^+}(X_1^*,X_2^*)
&=& 2(k+n) \Tr C_2C_1=(k+n)\, g^+_{p_0}(X_1^*,Y_1^*);\\
Ric_{p_0}^{g^+}(X_1^*,Y_1^*)&=&0  \quad \textrm{and} \quad  g^+_{p_0}(X_1^*,Y_1^*)=0;\\
Ric_{p_0}^{g^+}(Y_1^*,Y_2^*)
&=& \frac{n}{2}\Tr (B'_2B_1+B'_1B_2)
+\frac{1}{4}(4k\Tr(B'_1B_2+B'_2B_1)+2n\Tr(B_2B'_1+B_1B'_2)) \\
&=&(k+n)(\Tr B'_1B_2 + \Tr B'_2B_1) = (k+n)\,g^+_{p_0}(Y_1^*,Y_1^*).
\end{eqnarray*}}}
The almost complex structure $\J^-$ is special : $ChernRicci^{\J^-}=2(n-k)\, \omega$. Indeed, using the identities \eqref{eq:RicChern1}, \eqref{eq:RicChern2} and  \eqref{eq:RicChern3}, one has :
{\scriptsize{\begin{eqnarray*}
ChernRicci^{\J^-}_{p_0}({X_i^*},X_j^*)
     && 2n\Tr J_{2n}[C_i,C_j]-2k\Tr J_{2n}[C_i,C_j]=2(n-k)\,\omega_{p_0}({X_i^*},X_j^*);\\
ChernRicci^{\J^-}_{p_0}({X_i^*},Y_k^*)&=&0 \quad \textrm{and} \quad  \omega_{p_0}(X_1^*,Y_1^*)=0;\\
ChernRicci^{\J^-}_{p_0}({Y_i^*},Y_j^*)&=& (2n-2k)\Tr J_{2n}(B'_iB_j-B'_jB_i)=2(n-k)\,\omega_{p_0}({Y_i^*},Y_j^*).
\end{eqnarray*}}}

%%%%%%%%%%%%%%%%%%%%%%%%%%%%%%%%%%%%%%%

\subsection{$G=O(k+2n,\R),\, O(k,2n,\R),\,SO(k+2n,\R)$ or $ SO_0(k,2n,\R)$ with $k,n\geq1$}

The orthogonal groups $O(k+2n,\R)=\left\{ g\in Gl(k+2n,\R)\, \vert \, ^{tr}gg=\id\right\}$, the pseudo-orthogonal  groups $O(k,2n,\R)=\left\{ g\in Gl(k+2n,\R)\, \vert \, ^{tr}g{\tiny{\begin{pmatrix}
\Id_k & 0\cr
0 &  -\Id_{2n}
\end{pmatrix}}} g = {\tiny{\begin{pmatrix}
\Id_k & 0\cr
0 &  -\Id_{2n}
\end{pmatrix}}} \right\}$, 
and their connected components  are $\tilde{\sigma}$-stable, since $R\in O(k,2n,\R)\cap O(k+2n,\R)$.
Clearly, $\rho$ belongs to the Lie algebras $\mathfrak{o}(k+2n,\R)=\{ X\in \mathfrak{gl}(k+2n,\R)\, \vert \, ^{tr}X+X=0\} $ and $\mathfrak{o}(k,2n,\R)=\{ X\in \mathfrak{gl}(k+2n,\R)\, \vert \, ^{tr}X{\tiny{\begin{pmatrix}
\Id_k & 0\cr
0 &  -\Id_{2n}
\end{pmatrix}}}+{\tiny{\begin{pmatrix}
\Id_k & 0\cr
0 &  -\Id_{2n}
\end{pmatrix}}}X=0\} $.\\
  These algebras are simple for $k+2n\neq 4$ and in each case the 
 map $\beta: \g \times \g \rightarrow \R $ is non degenerate. The subspace $\p=\g^\sigma_{-1}=\left\{ \left(\begin{smallmatrix}
0&B\\ B'& 0
\end{smallmatrix}\right) \in \g\right\}$ does not vanish. Hence the $4$-symmetric space $G/G^{\tilde{\sigma}}$ identifies with the coadjoint orbit of $\rho^\flat$ in the corresponding $\g^*$
 and the $2$-form $\omega$ is symplectic.
 Thus
\begin{align*}
&O(k+2n,\R)/O(k,\R)\times U(n)
\qquad SO(k+2n,\R)/SO(k,\R)\times U(n)\cr
&O(k,2n,\R)/O(k,\R)\times U(n) \qquad 
SO_0(k,2n,\R)/SO(k,\R)\times U(n),
\end{align*}
where $U(n)=\{ A\in O(2n,\R)\, \vert \, A J_{2n}=J_{2n} A\}$, endowed with $\omega$ and $\J^+$ are pseudo-K\"ahler $4$-symmetric spaces. Endowed with $\omega$ and $\J^-$, these are 
almost pseudo-K\"ahler $4$-symmetric spaces, with the invariant natural almost complex structure $\J^-$ being maximally non integrable.
Indeed, $[\p,\p]\cap\g^\sigma_{-1}=\g^\sigma_{-1}$ and $[\p,\g^\sigma_{-1}]= \p$,
using remarks \ref{rmkintegr1} and  \ref{rmkintegr},  and the fact that 
the $2n\times 2n$ matrices $C_{ij}= E_{i\, j}-E_{j\, i}-E_{n+i\, n+j}+E_{n+j\, n+i}$, for $1\le i<j\le n$,
(where $ E_{i\,j}$ is the $2n\times 2n$ matrix whose only non vanishing entry is a $1$ at the intersection of the $i^{th}$ row and the $j^{th}$ column), satisfy $\sum_{ij}C_{ij}^2=-(n-1)\Id_{2n}$, and ${\tiny{\begin{pmatrix}
0&0\\0&C_{ij}
\end{pmatrix}}}\in\mathfrak{g}^{\sigma}_{-1}$ since 
\begin{equation}\label{eq:g1sigma}
\mathfrak{g}^{\sigma}_{-1}=\left\{ \left(\begin{matrix}
0&0\\0&C
\end{matrix}\right)\,\vert\, C=\left(\begin{smallmatrix}  c & c'\cr c'&  -c\end{smallmatrix}\right),\,\, ^{tr}c=-c\quad ^{tr}c'=-c'\right\}.
\end{equation}  
For $\g=\mathfrak{o}(k+2n,\R)$, then $\mathfrak{g}^{\sigma}_{-1}+\p=
\left(\begin{smallmatrix} 0&  B\cr  -^{tr}B   &  C \end{smallmatrix}\right)$ with $B$ a $k\times 2n$-matrix and $C$ as in formula \ref{eq:g1sigma}. For such an $X \in \mathfrak{g}^{\sigma}_{-1}+\p$,
\begin{eqnarray*}
-\omega_{p_0}(\pi_*X,{\jj^\pm}_{p_0}\pi_*X)&=&-\tilde{\Omega}(X,J^\pm X)=2\Tr (^{tr}B B)\mp2 \Tr(C^2)\cr
&=&2\sum_{i=1}^k\sum_{j=1}^{2n} (B^i_j)^2\mp4 \Tr ({c}^2+{c'}^2)\cr
&=& 2\sum_{i=1}^k\sum_{j=1}^{2n}(B^i_j)^2 \pm4\sum_{i=1}^n\sum_{j=1}^{n}(({c}^i_j)^2+({c'}^i_j)^2).
\end{eqnarray*}
Thus, the adjoint orbit of $\rho$ in $\mathfrak{o}(k+2n,\R)$, $ SO(k+2n,\R)/SO(k,\R)\times U(n)$, endowed with $-\omega$ and the invariant complex structure $\J^+$, is a K\"ahler manifold.\\
The case where $k=1$ corresponds to the twistor space $SO(1+2n,\R)/ U(n)$ on the sphere $S^{2n}=SO(1+2n,\R)/SO(2n,\R)$. \\

For $\g=\mathfrak{o}(k,2n,\R)$, then $\mathfrak{g}^{\sigma}_{-1}+\p=\{
X=\left(\begin{smallmatrix} 0&  B\cr  ^{tr}B   &  C \end{smallmatrix}\right)\}$ with $B$ a $k\times 2n$-matrix and $C$ as in formula \ref{eq:g1sigma}. For such an $X \in \mathfrak{g}^{\sigma}_{-1}+\p$,
\begin{eqnarray*}
\omega_{p_0}(\pi_*X,{\jj^\pm}_{p_0}\pi_*X)&=&\tilde{\Omega}(X,J^\pm X)=2\Tr (^{tr}B B)\pm 2 \Tr(C^2)\cr
&=&2\sum_{i=1}^k\sum_{j=1}^{2n} (B^i_j)^2\pm 4 \Tr ({c}^2+{c'}^2)\cr
&=& 2\sum_{i=1}^k\sum_{j=1}^{2n}(B^i_j)^2 \mp 4\sum_{i=1}^n\sum_{j=1}^{n}(({c}^i_j)^2+({c'}^i_j)^2).
\end{eqnarray*}
Thus, on the adjoint orbit of $\rho$ in $\mathfrak{o}(k,2n,\R)$, $ SO_0(k,2n,\R)/SO(k,\R)\times U(n)$, endowed with $\omega$, the invariant structure $\J^-$, is positive admissible and maximally non integrable; we have a  maximally non integrable almost K\"ahler manifold.\\
The case where $k=1$ corresponds to the twistor space $SO_o(1,2n,\R)/ U(n)$ on the hyperbolic space  $SO_o(1,2n,\R)/SO(2n,\R)$. \\

In all cases, a basis of the  $C$'s corresponding to elements in $\mathfrak{g}^{\sigma}_{-1}$ is given by 
$$\{E_{i\, j}-E_{j\, i}-E_{n+i\, n+j}+E_{n+j\, n+i}, E_{n+i \, j}-E_{n+j \, i}+E_{i\, n+j}-E_{j\, n+i}, 1\le i<j\le n\}$$ so that, for any $D=\left(\begin{smallmatrix}
d & d'\cr
-d' & d
\end{smallmatrix}\right)$, and for any $C_1,C_2$ of the form $\left(\begin{smallmatrix}  c_i & c_i'\cr c_i'&  -c_i\end{smallmatrix}\right)$
{\scriptsize{ 
\begin{eqnarray*}
\Tr\Big( C\rightarrow \left((D)C+C (^{tr}D)\right)\Big)&=&\sum_{1\le i<j\le n}(d_{ii}+d_{jj}+d_{ii}+d_{jj})=2(n-1)\Tr d =(n-1)\Tr D\\
\Tr\Big( C\rightarrow \left( C_1C C_2+(^{tr}C_2)C (^{tr}C_1)\right)\Big)&=& 0.\\
\end{eqnarray*}
}}
Since $\p=\{ \left(\begin{smallmatrix}
0 & B\cr
B'=\mp^{tr}B &  0
\end{smallmatrix}\right)\}$ (with the $-$ sign for $\g=\mathfrak{o}(k+2n,\R)$ and the $+$ sign for $\g=\mathfrak{o}(k,2n,\R)$), a basis of the  set of pairs $(B,B')$ is given by 
$$\{(E_{r\, \ell}, \mp E_{\ell\,  r}), 1\le \ell \le 2n , 1\le r\le k \}$$ 
and, for any $G\in Mat(k\times k),\beta \beta'\in Mat(k\times 2n)$ and $F \in Mat(2n\times 2n)$,
{\scriptsize{
\begin{eqnarray*}
\Tr\Big((B,B'=\mp \,^{tr}B)&\longrightarrow &(G B F +\beta B'\beta', \mp\,^{tr}( G B F +\beta B'\beta')\Big)\\
&=&\sum_ {  1\le r \le k,1\le \ell \le 2n}(G_{rr}F_{\ell\ell}\mp \beta_{r\ell}\beta'_{r\ell})=\Tr G\Tr F \mp \Tr(\,^{tr}\beta \beta').
\end{eqnarray*}}}

The metric $g^+$ is Einstein $Ric^{g^+}=\half(k+n-1)g^+$.Indeed, the identities \eqref{eq:Einstein1}, \eqref{eq:Einstein2} and  \eqref{eq:Einstein3} become :
{\scriptsize{ \begin{eqnarray*}
Ric_{p_0}^{g^+}(X_1^*,X_2^*)
&=&\Tr\left( C\rightarrow (C_2C_1C+CC_1C_2-C_1CC_2-C_2CC_1)\right)+\Tr \left( (B,B')\rightarrow (BC_1C_2, C_2C_1B')\right)\\
&=&(n-1)\Tr C_2C_1+ k \Tr C_2C_1=\half(k+n-1)g^+_{p_0}(X_1^*,X_2^*)\\
Ric_{p_0}^{g^+}(X_1^*,Y_1^*)&=&0  \quad \textrm{and} \quad  g^+_{p_0}(X_1^*,Y_1^*)=0;\\
Ric_{p_0}^{g^+}(Y_1^*,Y_2^*)&=&\half\Tr\left( C\rightarrow \half((B'_2B_1-J_{2n}B'_2B_1J_{2n})C+C(B'_1B_2-J_{2n}B'_1B_2J_{2n})\right)\\
&&+\frac{1}{4}\Tr\Bigg( (B,B')\rightarrow \Big( (2B_2B'_1-B_1B'_2)B+B(3B'_1B_2+B'_2B_1-JB'_1B_2J +JB'_2B_1J)\nonumber\\
&&\qquad\qquad \qquad  -(B_1JB'_2+2B_2JB'_1)BJ-2B_2B'B_1-3B_1B'B_2+B_1JB'B_2J+2 B_2JB'B_1J  \, ,\nonumber\\
&&\qquad\qquad \qquad  B'(2B_1B'_2-B_2B'_1)+(3B'_2B_1+B'_1B_2+JB'_1B_2J -JB'_2B_1J) B'      \nonumber   \\
&&\qquad\qquad \qquad  -JB'(B_2JB'_1+2B_1JB'_2)-2B'_1BB'_2-3B'_2BB'_1+JB'_2BJB'_1+2 JB'_1BJB'_2 \Big)\Bigg)\\
&=& \half (n-1) \Tr(\half((B'_2B_1-J_{2n}B'_2B_1J_{2n}))\\
&&+\frac{1}{4}(2n\Tr(2B_2B'_1-B_1B'_2)+k\Tr(3B'_1B_2+B'_2B_1-JB'_1B_2J +JB'_2B_1J))\\
&&+\frac{1}{4}\Tr(-2B'_2B_1-3B'_1B_2-JB'_1B_2J-2JB'_2B_1J)\\
&=&  (\half (n-1) +\frac{1}{4}(2n+ 4k-2))   \Tr (B'2B_1)= (n-1+k)\half g^+_{p_0}(Y_1^*,Y_2^*).
 \end{eqnarray*}
 }}

The almost complex structure $\J^-$ is special : $ChernRicci^{\J^-}=(n-1-k)\omega$. This results directly from the identities \eqref{eq:RicChern1}, \eqref{eq:RicChern2} and  \eqref{eq:RicChern3}, the computations given above, and the fact that $D(X_i,X_j)=J_{2n}[C_i,C_j]$ and 
$D(Y_i,Y_j) =\half( J_{2n}(B'_iB_j-B'_jB_i)+(B'_iB_j-B'_jB_i)J_{2n})$ are both symmetric.
%%%%%%%%%%%%%%%%%%%%%%%%%%

%%%%%%%%%%%%%%%%%%%%%%%%%%%%%%%%

\subsection{  $G=U(k',n)$ or $U(k'+n)$}

The groups $U(k',n)$ (resp $U(k'+n)$) are seen as the subgroups  of the elements of $O(2k'+2n,\R)$ (resp. $O(2k'+2n,\R)$) which commute with $\left(\begin{smallmatrix} J_{2k'}&  0\cr  0   &  J_{2n} \end{smallmatrix}\right)$.

The element $\rho$ belongs to their Lie algebra and 
$\beta: \g \times \g \rightarrow \R : X,Y \mapsto \Tr(XY)$ is non degenerate. In both cases, 
$\g^\sigma_{-1}=\{ 0 \}$ and $G^{\tilde\sigma}=G^{\tilde\sigma^2}$. The $4$-symmetric spaces are thus  the symmetric spaces, 
$$
U(k'+n)/\left(U(k')\times U(n)\right)\qquad  U(k',n)/\left(U(k')\times U(n)\right)
$$ 
which identify with the coadjoint orbit of the element $\rho^\flat: \g \rightarrow \R : X \mapsto \Tr \rho X$. There is only one natural almost complex structure (up to sign), $\J^+$, defined by $\sigma$ on $\p$; it is integrable. 

For $\g=\mathfrak{u}(k'+n)$, $\p=\{\left(\begin{smallmatrix} 0&  B\cr  -^{tr}B  &  0 \end{smallmatrix}\right)\,\vert\, B= \left(\begin{smallmatrix} b&  b'\cr  -b'  &  b \end{smallmatrix}\right)\}$, and for such a $X\in \p$: 
\begin{eqnarray*}
-\omega_{p_0}(\pi_*X,{\jj^+}_{p_0}\pi_*X)&=&-\tilde{\Omega}(X,J^+ X)=-2\Tr (-^{tr}B B)\cr
&=& 2\sum_{i=1}^k\sum_{j=1}^{2n}(B^i_j)^2 
\end{eqnarray*}
so that $(-\omega, \J^+)$ define a K\"ahler structure on  $U(k'+n)/\left(U(k')\times U(n)\right)$.

For $\g=\mathfrak{u}(k',n)$, $\p=\{\left(\begin{smallmatrix} 0&  B\cr  ^{tr}B  &  0 \end{smallmatrix}\right)\,\vert\, B= \left(\begin{smallmatrix} b&  b'\cr  -b'  &  b \end{smallmatrix}\right)\}$, and for such a $X\in \p$: 
\begin{eqnarray*}
\omega_{p_0}(\pi_*X,{\jj^+}_{p_0}\pi_*X)&=&\tilde{\Omega}(X,J^+ X)=2\Tr (^{tr}B B)\cr
&=& 2\sum_{i=1}^k\sum_{j=1}^{2n}(B^i_j)^2 
\end{eqnarray*}
so that $(\omega, \J^+)$ define a K\"ahler structure on  $U(k',n)/\left(U(k')\times U(n)\right)$.

In both case, the action of the stabilizer $U(k')\times U(n)$ on $\p\simeq Mat(k'\times n,\C)\simeq \C^{k'}\otimes\C^n$ is irreducible so the K\"ahler metric is Einstein.

The curvature at the base point is given by $R_{p_0}^{g}({Y_1^*},Y_2^*)Y_3^*=-([[Y_1,Y_2],Y_3])^*_{p_0}$, so  
{\scriptsize{
\begin{eqnarray*}
Ric_{p_0}^{g}(Y_1^*,Y_2^*)&=&\Tr (\ad Y_2\circ \ad Y_1)\vert_{\p} \\
&=&\Tr\vert_{\p}\Big( B\rightarrow (B_2B'_1B+BB'_1B_2-B_2B'B_1-B_1B'B_2)\Big)\,\,\textrm{ with }\,\, B'=\mp \,^{tr}B\\
&=& \Tr(b,b')\rightarrow \Big((b_2\tilde b_1-b'_2\tilde b'_1)b+b(\tilde b_1b_2-\tilde b'_1b'_2)-b_2\tilde b b_1-b_1\tilde b b_2+b'_2\tilde b b'_1+b'_1\tilde b b'_2 \, ,\\
   &&\qquad \qquad \quad(b_2\tilde b_1-b'_2\tilde b'_1)b'+b'(\tilde b_1b_2-\tilde b'_1b'_2)+b_2\tilde b' b_1+b_1\tilde b' b_2-b'_2\tilde b' b'_1-b'_1\tilde b' b'_2
   \Big)\\
&=&(n+k') \Tr(B_2B'_1)= \half (n+k') \Tr(B'_1B_2+B'_2B_1).
\end{eqnarray*}
avec $b,b'\in Mat(k'\times n)$ and $B'=\left(\begin{smallmatrix} \tilde b&  \tilde b'\cr  -\tilde b'  &  \tilde b \end{smallmatrix}\right)$.
}}
%%%%%%%%%%%%%%%%%%%%%%%%%%
\subsection{$G=Sp(k+2n,\R)\simeq Sp(\mathbb{R}^{2(k'+n)},\widetilde\Omega)$}

We set $k=2k'$ and we  choose a basis of $\R^{k+2n}=\R^{2(k'+n)}$ in which the non-degenerate  skew-symmetric  $2$-form writes
 $$
 \widetilde\Omega= \begin{pmatrix}
\Omega_{2k'}&0\\0&\Omega_{2n}
\end{pmatrix}\qquad \textrm{ with} \quad \Omega_{2r}:=-J_{2r}= \begin{pmatrix}
O&\Id_r\\ -\Id_r  & 0
\end{pmatrix}.
$$
The subgroup $Sp(\R^{2(k'+n)},\widetilde\Omega)=\left\{ g\in Gl(k+2n,\R)\, \vert \,^{tr}g\,\widetilde\Omega \, g =\widetilde\Omega\right\}$ is $\tilde{\sigma}$-stable, because $R\in Sp(\R^{2(k'+n)},\widetilde\Omega)$. The element $\rho$ is in the Lie algebra $\mathfrak{sp}(\R^{2(k'+n)},\widetilde\Omega)$:
$$\mathfrak{sp}(\R^{2(k'+n)},\widetilde\Omega)=\left\{\left(\begin{smallmatrix}
A\,\,&\,\Omega_{2k'} \,^{tr}B \Omega_{2n}\cr B\,\,&C
\end{smallmatrix}\right)\,\vert\,  ^{tr}A\Omega_{2k'}+\Omega_{2k'} A=0,\,^{tr}C\Omega_{2n}+\Omega_{2n} C=0         \right\}.$$
This Lie algebra $\mathfrak{sp}(\R^{2(k'+n)},\widetilde\Omega)\simeq\mathfrak{sp}(k+2n,\R)$ is simple, thus the $2$-form $\beta$  is non degenerate. The subspace $\p$ is given by
$$ \p=\left\{\begin{pmatrix}
0\,\,&\,B\cr 
\Omega_{2n} \,^{tr}B \Omega_{2k'}\,\,&0
\end{pmatrix}\,\vert\,  B\in \textrm{Mat}(2n\times 2k')\, \right\}\neq \lbrace 0 \rbrace .$$
The $4$-symmetric space
$$
G/G_{0}^{\tilde{\sigma}}=Sp(\R^{2(k'+n)},\widetilde\Omega)/(Sp(\R^{2k'},\Omega_k)\times U(n)),
$$ 
%o\`u l'on voit le groupe unitaire $U(n)$ comme l'ensemble des matrices de $Sp(\R^{2n},\Omega_{2n})$ qui commutent avec $J_{2n}$;
is the adjoint orbit of $\rho$ in $\mathfrak{sp}(\R^{2(k'+n)},\widetilde\Omega)$, or the coadjoint orbit of $\rho^\flat$ in the dual Lie algebra, and the form $\omega$ coincides with the canonical symplectic form of Kirillov-Kostant-Souriau. 

The almost complex structure $\jj^+$ is an invariant complex structure and the almost complex structure $\jj^-$ is maximally non integrable since $\g$ is simple and 
\begin{eqnarray*}
\mathfrak{g}^{\sigma}_{-1}&=&\left\{\begin{pmatrix}
0&0\cr 0&C
\end{pmatrix}\,\vert\,  ^{tr}C\Omega_{2n}+\Omega_{2n} C=0  ~\textrm{ and }~ CJ_{2n}+J_{2n}C=0 \right\},\\
&=&\left\{\begin{pmatrix}
0&0\cr 0&C\end{pmatrix}\,\vert\, C= \left(\begin{smallmatrix}
c& c'\cr c'&-c\end{smallmatrix}\right) ~~ \textrm{with}~~ ^{tr}c=c~\textrm{ and }~ ^{tr}c'=c'\, \right\}
\end{eqnarray*}
contains the element  $\left(\begin{smallmatrix}
0&0\\0&C
\end{smallmatrix}\right) $ where $C= \left(\begin{smallmatrix}
\Id_n&0\\0&-\Id_n
\end{smallmatrix}\right)$ satisfies 
$C^2=\Id_{2n}$.

A basis of the  $C$'s corresponding to elements in $\mathfrak{g}^{\sigma}_{-1}$ is given by 
$$\{E_{i\, j}+E_{j\, i}-E_{n+i\, n+j}-E_{n+j\, n+i}, E_{n+i \, j}+E_{n+j \, i}+E_{i\, n+j}+E_{j\, n+i}, 1\le i\le j\le n\}$$ so that, for any $D=\left(\begin{smallmatrix}
d & d'\cr
-d' & d
\end{smallmatrix}\right)$,  
$$\Tr\left( C\rightarrow DC+CD\right)=\sum_{1\le i \le j\le n}(d_{ii}+d_{jj}+d_{ii}+d_{jj})=2(n+1)\Tr d=(n+1)\Tr D.$$
Since $\p=\{ \left(\begin{smallmatrix}
0 & B\cr
B'=\Omega_{2n} \,^{tr}B \Omega_{2k'}&  0
\end{smallmatrix}\right)\}$, a basis of the  $(B,B')$ is given by 
$$\{(E_{r\, \ell}, \Omega_{2n}E_{\ell\,  r}\Omega_{2k'}), 1\le \ell \le 2n , 1\le r\le 2k' \} $$and for any $D$ as above such that $^{tr}D\,\Omega_{2n}=\Omega_{2n}D$,
$$\Tr\left((B,B')\rightarrow (BD, DB')\right)=\sum_ {1\le \ell \le 2n , 1\le r \le 2k'}
D_{\ell \ell}=2k' \Tr D.$$
The almost complex structure $\J^-$ is special : $ChernRicci^{\J^-}=(n+1-k)\omega$.

The metric $g^+$ is Einstein $Ric^{g^+}=\half(2k'+n+1)g^+$. Indeed, the identities \eqref{eq:Einstein1}, \eqref{eq:Einstein2} and  \eqref{eq:Einstein3} become :
{\scriptsize{ \begin{eqnarray*}
Ric_{p_0}^{g^+}(X_1^*,X_2^*)
&=&\Tr\left( C\rightarrow (C_2C_1C+CC_1C_2-C_1CC_2-C_2CC_1)\right)+\Tr \left( (B,B')\rightarrow (BC_1C_2, C_2C_1B')\right)\\
&=&(n+1)\Tr C_2C_1+ 2k' \Tr C_2C_1=\half(2k'+n+1)g^+_{p_0}(X_1^*,X_2^*)\\
Ric_{p_0}^{g^+}(X_1^*,Y_1^*)&=&0  \quad \textrm{and} \quad  g^+_{p_0}(X_1^*,Y_1^*)=0;\\
Ric_{p_0}^{g^+}(Y_1^*,Y_2^*)&=&\half\Tr\left( C\rightarrow \half((B'_2B_1-J_{2n}B'_2B_1J_{2n})C+C(B'_1B_2-J_{2n}B'_1B_2J_{2n})\right)\\
&&+\frac{1}{4}\Tr\Bigg( B\rightarrow \Big( (2B_2B'_1-B_1B'_2)B+B(3B'_1B_2+B'_2B_1-JB'_1B_2J +JB'_2B_1J)\nonumber\\
&&\qquad\qquad \qquad  -(B_1JB'_2+2B_2JB'_1)BJ-2B_2B'B_1-3B_1B'B_2+B_1JB'B_2J+2 B_2JB'B_1J \Big)   \Bigg)\\
&=& \half (n+1) \Tr(\half((B'_2B_1-J_{2n}B'_2B_1J_{2n}))\\
&&+\frac{1}{4}(2n\Tr(2B_2B'_1-B_1B'_2)+2k'\Tr(3B'_1B_2+B'_2B_1-JB'_1B_2J +JB'_2B_1J))\\
&&+\frac{1}{4}\Tr(2\,^{tr}B_1J_{2k}B_2J_{2n}+3\,^{tr}B_2JB_1J-J\,^{tr}B_2JB_1-2J\,^{tr}B_1JB_2)
\\
&=&  (\half (n+1) +\frac{1}{4}(2n+ 8k'+2))   \Tr (B'2B_1)= (n+1+2k')\half g^+_{p_0}(Y_1^*,Y_2^*).
 \end{eqnarray*}
 }}

\end{document}